\input amstex\documentstyle {amsppt}  
\pagewidth{12.5 cm}\pageheight{19 cm}\magnification\magstep1
\topmatter
\title Character sheaves on disconnected groups, I\endtitle
\author G. Lusztig\endauthor
\address Department of Mathematics, M.I.T., Cambridge, MA 02139\endaddress
\thanks Supported in part by the National Science Foundation\endthanks
\endtopmatter   
\document

\define\Lie{\text{\rm Lie }}

\define\frl{\forall}

\define\si{\sim}

\define\sqc{\sqcup}

\define\qua{\quad}

\define\hL{\hat L}

\define\bC{\bar C}
\define\bY{\bar Y}

\define\op{\oplus}

\define\em{\emptyset}
\define\imp{\implies}
\define\ra{\rangle}

\define\iy{\infty}
\define\m{\mapsto}
\define\do{\dots}
\define\la{\langle}
\define\bsl{\backslash}

\define\lra{\leftrightarrow}

\define\sub{\subset}
\define\bxt{\boxtimes}
\define\T{\times}
\define\ti{\tilde}
\define\nl{\newline}
\redefine\i{^{-1}}
\define\fra{\frac}
\define\un{\underline}

\define\ot{\otimes}
\define\bbq{\bar{\QQ}_l}

\define\Ad{\text{\rm Ad}}
\define\Hom{\text{\rm Hom}}

\define\Aut{\text{\rm Aut}}

\define\Ker{\text{\rm Ker}}

\define\supp{\text{\rm supp}}

\define\a{\alpha}
\redefine\b{\beta}
\redefine\c{\chi}
\define\g{\gamma}
\redefine\d{\delta}
\define\e{\epsilon}

\redefine\o{\omega}
\define\p{\pi}
\define\ph{\phi}

\define\r{\rho}
\define\s{\sigma}
\redefine\t{\tau}

\redefine\l{\lambda}
\define\z{\zeta}

\redefine\G{\Gamma}

\define\boc{\bold c}

\define\kk{\bold k}

\define\zz{\bold z}

\redefine\AA{\bold A}

\define\QQ{\bold Q}

\define\ZZ{\bold Z}

\define\cc{\Cal C}
\define\cd{\Cal D}
\define\ce{\Cal E}

\define\ch{\Cal H}

\define\cn{\Cal N}
\define\co{\Cal O}
\define\cp{\Cal P}

\define\cs{\Cal S}

\define\cv{\Cal V}
\define\cw{\Cal W}
\define\cz{\Cal Z}
\define\cx{\Cal X}

\define\fg{\frak g}

\define\tp{\ti p}

\define\tG{\ti G}
\define\tH{\ti H}

\define\tL{\ti L}

\define\tP{\ti P}

\define\tT{\ti T}

\define\tY{\ti Y}

\define\sps{\supset}

\define\app{\asymp}
\define\utP{\un{\tP}}
\define\uP{\un P}
\define\bnu{\bar\nu}
\define\bc{\bar c}

\define\br{\bar r}
\define\bg{\bar g}

\define\bE{\bar E}
\define\bS{\bar S}
\define\bP{\bar P}
\define\bce{\bar\ce}
\define\tce{\ti\ce}
\define\BO{B}
\define\DS{D}
\define\FU{L1}
\define\IC{L2}
\define\CS{L3}
\define\IN{L4}
\define\CL{L5}
\define\SP{Sp}
\define\ST{St}

\head Introduction\endhead
Our aim in this series of papers is to develop a theory of character sheaves on a not 
necessarily connected reductive algebraic group $G$. In the case of connected groups 
such a theory appeared in \cite{\IC}, \cite{\CS}. An extension to disconnected 
groups has been sketched in \cite{\IN} without proofs; here we try to give a fuller
and more precise treatment and to supply the proofs that were missing in \cite{\IN}.
The main object of the theory, the character sheaves of $G$, are certain simple 
perverse sheaves on $G$, equivariant with respect to the conjugation action of the 
identity component of $G$. At least for connected $G$ (over a finite field) the 
character sheaves are intimately related with the characters of irreducible 
representations of the group of rational points of $G$ and such a relationship is 
also expected in the disconnected case. The theory of character sheaves on $G$ is 
also crucial for the classification of "unipotent representations" of simple $p$-adic 
groups \cite{\CL} and here one is forced to allow $G$ to be disconnected if one wants 
to include $p$-adic groups that are not inner forms of split groups. 

The present paper tries to extend parts of \cite{\IC,\S1-\S4} from the connected case
to the general case. We develop enough background so that we are able to define (see 
6.7) the notion of "admissible complex" on $G$, one of the two incarnations of the
character sheaves of $G$. One of the themes of this paper is the construction of a 
decomposition of $G$ into finitely many strata (generalizing a construction in 
\cite{\IC, 3.1}), see \S3. Each stratum is a locally closed, irreducible, smooth 
subvariety of $G$. Each stratum is a union of $G^0$-conjugacy classes of fixed 
dimension; more precisely, the centralizers of two points in the same stratum have 
$G^0$-conjugate identity components. Also, the closure of any stratum is a union of 
strata. Each stratum carries some natural local systems which extend to intersection 
cohomology complexes on the closure, which we also describe by a direct image 
construction, using the dimension estimates in \S4. 

I wish to thank F. Digne and J. Michel for some useful comments on \cite{\IN}.

\head Contents\endhead
1. Preliminaries on reductive groups.

2. Isolated elements of $G$.

3. A stratification of $G$.

4. Dimension estimates.

5. Some complexes on $G$.

6. Cuspidal local systems.

\head 1. Preliminaries on reductive groups\endhead
\subhead 1.1\endsubhead
We fix an algebraically closed field $\kk$. All algebraic variety are assumed to be
over $\kk$. All algebraic groups are assumed to be affine. 

We shall use the following notation. If $H$ is a group, the centre of $H$ is denoted by
$\cz_H$; if $H'$ is a subgroup of $H$, let $N_H(H')=\{h\in H;hH'h\i=H'\}$. If in 
addition $H''$ is a subgroup of $H$, let 
$Z_{H'}(H'')=\{h'\in H';h'h''=h''h'\qua \frl h''\in H''\}$; if $h\in H$, let
$Z_{H'}(h)=\{h'\in H';h'h=hh'\}$. 
If $H$ is an algebraic group, we denote by $H^0$ the identity component of $H$ and we set
$H_{ss}=H/\cz_{H^0}^0$; for $h\in H$ we denote by $h_s$ (resp. $h_u$) the semisimple 
(resp. unipotent) part of $h$, so that $h=h_sh_u=h_uh_s$. If $X$ is a subset of $H$ we 
set $X_s=\{h_s;h\in X\}$. The unipotent radical of $H$ (assumed to be connected) is 
denoted by $U_H$. 

We fix an algebraic group $G$ such that $G^0$ is reductive. (We then say that $G$ is 
reductive.) Let $\fg=\Lie G$.

\subhead 1.2\endsubhead 
Let $T$ be a torus and let $f:T@>>>T$ be an automorphism of finite order with fixed 
point set $T^f$.We show that

(a) {\it the homomorphism $(T^f)^0\T T@>>>T,(t,x)\m xtf(x)\i$ is surjective.}
\nl
This can be reduced to an analogous statement about a finite dimensional $\QQ$-vector 
space $V$ and a linear map $\ph:V@>>>V$ of finite order: the linear map 
$\Ker(\ph-1)\T V@>>>V,(w,v)\m w+v-\ph(v)$ is surjective. 

\subhead 1.3\endsubhead
Let $g\in G$. We show that

(a) {\it the homomorphism 
$(\cz_{G^0}^0\cap Z_G(g))^0\T\cz_{G^0}^0@>>>\cz_{G^0}^0,(t,x)\m xtgx\i g\i$
is surjective}. 
\nl
$\Ad(g):\cz_{G^0}^0@>>>\cz_{G^0}^0$ is of finite order since some power of
$g$ is in $G^0$. Therefore (a) is a special case of 1.2(a).

\subhead 1.4\endsubhead
Let $g\in G$. Then $g$ normalizes some Borel of $G^0$, see \cite{\ST, 7.2}. Following 
\cite{\ST, 9}, we say that $g$ is {\it quasi-semisimple} if there exist a Borel $B$ of 
$G^0$ and a maximal torus $T$ of $B$ such that $gBg\i=B,gTg\i=T$. If $g$ is semisimple then 
it is quasi-semisimple, see \cite{\ST, 7.5, 7.6}. More generally, by an argument 
similar to that in \cite{\ST, 7.6}, we see that

(a) {\it if $g$ is semisimple and $P$ is a parabolic of $G^0$ such that $gPg\i=P$ then
there exists a Levi $L$ of $P$ such that} $gLg\i=L$.
\nl
Here are some further results.

(b) {\it If $g$ is semisimple or, more generally, quasi-semisimple then $Z_G(g)$ is 
reductive.} See \cite{\SP, 1.17, 2.21}.

(c) $g$ {\it is quasi-semisimple if and only if $g_u$ is quasi-semisimple in the 
reductive group} $Z_G(g_s)$. See \cite{\SP, 2.22}.

(d) {\it If $g$ is quasi-semisimple and $T_1$ is a maximal torus of $Z_G(g)^0$ then 
there is a unique maximal torus $T$ of $G^0$ such that $T_1\sub T$.  See 
\cite{\SP, 1.15}. (We have necessarily} $T=Z_{G^0}(T_1),T_1=(T\cap Z_G(g))^0$.)

(e) {\it $g$ is quasi-semisimple if and only if the $G^0$-conjugacy class of $g$ is 
closed in $G$.} See \cite{\SP, 1.15}.

\subhead 1.5\endsubhead
Assume that $s\in G$ is semisimple. Let $T_1$ be a maximal torus of $Z_G(s)^0$. 
Clearly, some power of $s$ is in $\cz_{Z_G(s)^0}$ hence in $T_1$ (since $Z_G(s)^0$ is 
reductive). Thus, the subgroup $\la s\ra T_1$ generated by $s$ and $T_1$ is a closed 
diagonalizable subgroup of $G$ with identity component $T_1$. For any character 
$\a:\la s\ra T_1@>>>\kk^*$ let 
$\fg_\a=\{x\in\fg;\Ad(a)x=\a(a)x\quad \frl a\in\la s\ra T_1\}$. Then 
$\fg=\op_\a\fg_\a$. Let $R$ be the set of all $\a$ such that $\a\ne 1,\fg_\a\ne 0$. For
$\a\in R$ we have necessarily $\dim\fg_a=1$. (Compare \cite{\CL, 6.18}.) We have
$\fg_1=\Lie T$ where $T=Z_{G^0}(T_1)$.

\subhead 1.6\endsubhead
Let $\bE(G^0)$ be the set of all pairs $(B,T)$ where $B$ is a Borel of $G^0$ and $T$ is
a maximal torus of $B$. The group $\Aut(G^0)$ of automorphisms of $G^0$ acts naturally 
on $\bE(G^0)$. It is known that to $G^0$ one can associate canonically an algebraic 
variety $E(G^0)$ whose points are called "\'epinglages" with the following properties. 

(i) There is a natural action of $\Aut(G^0)$ on $E(G^0)$ which restricts to a free
transitive action of the group of inner automorphisms of $G^0$ on $E(G^0)$.

(ii) There is a natural $\Aut(G^0)$-equivariant map $p:E(G^0)@>>>\bE(G^0)$.

\subhead 1.7\endsubhead
Let $g\in G$. Assume that $\Ad(g)e=e$ where $e\in E(G^0)$. Let $(B,T)=p(e)$. Then 
$gBg\i=B,gTg\i=T$ hence $g$ is quasi-semisimple. In particular, $Z_G(g)$ is reductive. 
The following results are known.

(a) $Z_B(g)^0$ is a Borel of $Z_G(g)^0$.

(b) $\cz_{Z_G(g)^0}^0=(\cz_{G^0}\cap Z_G(g))^0$;

(c) $P\m Z_P(g)^0$ {\it is a bijection between the set of parabolics of $G^0$ that 
contain $B$ and are normalized by $g$ and the set of parabolics of $Z_G(g)^0$ that 
contain $Z_B(g)^0$; moreover, if $L$ is a Levi of $P$ then $Z_L(g)^0$ is a Levi of }
$Z_P(g)^0$.

\subhead 1.8\endsubhead
Assume that $u\in G$ is unipotent and $(B,T)\in\bE(G^0)$ is such that
$uBu\i=B,uTu\i=T$. We show that 

(a) {\it there exists $e\in E(G^0)$ such that} $p(e)=(B,T),\Ad(u)e=e$.
\nl
Let $E'=\{e\in E(G^0);p(e)=(B,T)\}$. Then $T$ acts transitively on $E'$. Moreover, 
$\Ad(u)E'=E'$. Hence the subgroup $\la u\ra T$ generated by $u$ and $T$ acts on $E'$. 
Since $u\in N_GT$ and $(N_GT)^0=T$, some power of $u$ belongs to $T$. Hence 
$\la u\ra T$ is a closed subgroup of $G$ with identity component $T$. Let $e_0\in E'$. 
We have $\Ad(u)e_0=\Ad(t\i)e_0$ for some $t\in T$. Thus, $tu$ 
belongs to the stabilizer of $e_0$ in $\la u\ra T$, a closed subgroup of $\la u\ra T$. 
Then the unipotent part $(tu)_u$ also belongs to this stabilizer that is, 
$\Ad((tu)_u)e_0=e_0$. The image of $(tu)_s$ in the unipotent group $(\la u\ra T)/T$ 
must be $1$. Hence $(tu)_u$ has the same image as $tu$ or as $u$. Thus, $(tu)_u=t'u$ 
where $t'\in T$ and we have $\Ad(t'u)e_0=e_0$ with $t'u$ unipotent. By 1.2(a) with
$f:T@>>>T,f(x)=uxu\i$ (of finite order), we have $t'=t_2t_1ut_1\i u\i$ for some 
$t_1,t_2\in T$ with $ut_2=t_2u$. Then $t'u=t_2t_1ut_1\i$. Since $t_2$ is semisimple and
it commutes with $t_1ut_1\i$ which is unipotent, we see that $t_1ut_1\i=(t'u)_u=t'u$ 
we see that $t'u=t_1ut_1\i$. Thus $\Ad(t_1ut_1\i)e_0=e_0$ hence 
$\Ad(u)e=e$ where $e=\Ad(t_1)\i e_0$. This proves (a). 

\subhead 1.9\endsubhead
Let $D$ be a connected component of $G$ which contains some unipotent element of $G$. Then

(a) $D$ {\it contains a unique closed unipotent $G^0$-conjugacy class; this is the 
set of unipotent, quasi-semisimple elements in $D$.} See \cite{\SP, 2.21}.

\subhead 1.10\endsubhead
Let $P$ be a parabolic of $G^0$ and let $L$ be a Levi of $P$. Let 
$g\in N_GL\cap\in N_GP$. We show that

(a) $Z_{G^0}((\cz_L\cap Z_L(g))^0)=L$.
\nl
Let $L'=Z_{G^0}((\cz_L\cap Z_L(g))^0)$. Then $L'$ is a reductive, connected subgroup of
$G^0$ and $L\sub L'$. Moreover, $P\cap L'$ is a parabolic subgroup of $L'$ with Levi 
$L$. We have $gL'g\i=L'$, $g(P\cap L')g\i=P\cap L'$. If (a) is true for
$N_GL',P\cap L',L,g$ instead of $G,L,P,g$, then we would have 
$Z_{L'}((\cz_L\cap Z_L(g))^0)=L$. Hence $L'=L$. Thus, to prove (a), we may assume that 
$L'=G^0$ and we must show that $L=G^0$. Replacing $g$ by a left $L$-translate, we may 
assume that there exists $(B_1,T)\in\bE(L)$ such that $gB_1g\i=B_1,gTg\i=T$. Since $g$ 
normalizes $U_P$, it also normalizes the Borel $B=B_1U_P$ of $G^0$.
Let $e\in E(G^0)$ be such that $p(e)=(B,T)$. Then $p(\Ad(g)e)=(B,T)$. We can
find $g_0\in G^0$ such that $\Ad(g)e=\Ad(g_0)e$ hence $(B,T)=(g_0Bg_0\i,g_0Tg_0\i)$ and
$g_0\in T$, $\Ad(g_0\i g)e=e$. Replacing $g$ by $g_0\i g$ we may assume in addition 
that $\Ad(g)e=e$. Then automatically $\Ad(g)$ fixes an \'epinglage of $L$ which lies 
over $(B_1,T)\in\bE(L)$. By 1.7(b) for $(g,G)$ and $(g,L)$, we have 
$$(\cz_{G^0}\cap Z_G(g))^0=\cz_{Z_{G^0}(g)^0}^0,(\cz_L\cap Z_L(g))^0=\cz_{Z_L(g)^0}^0.
$$
Since $L'=G^0$, we have $(\cz_L\cap Z_L(g))^0\sub(\cz_{G^0}\cap Z_G(g))^0$, hence
$\cz_{Z_L(g)^0}^0\sub\cz_{Z_{G^0}(g)^0}^0$. Since $Z_L(g)^0$ is a Levi of a parabolic 
of $Z_{G^0}(g)^0$ (see 1.7(c)), we have $Z_L(g)^0=Z_{Z_{G^0}(g)^0}(\cz_{Z_L(g)^0}^0)$ 
hence $Z_L(g)^0=Z_{G^0}(g)^0$. It follows that $Z_P(g)^0=Z_{G^0}(g)^0$. Using 1.7(c) we
deduce that $P=G^0$ hence $L=G^0$. This proves (a).

\subhead 1.11\endsubhead
Let $P$ be a parabolic of $G^0$. Let $s\in N_GP$ be semisimple. Then

(a) $U_P\cap Z_G(s)$ {\it is connected.}
\nl
Since $s$ normalizes $U_P$, this follows from \cite{\BO, 9.3}.

\subhead 1.12\endsubhead
Let $P$ be a parabolic of $G^0$. Let $s\in N_GP$ be semisimple. Choose a Levi $L$ of 
$P$ such that $s\in N_GL$. (See 1.4(a).) Let $Q=Z_G(s)^0\cap P$. We show:

(a) $Q$ {\it is a parabolic of $Z_G(s)^0$ with Levi $Z_L(s)^0$ and }
$U_Q=U_P\cap Z_G(s)$;

(b) $Z_L(s)^0$ {\it and $Z_G(s)^0$ have a common maximal torus and}
$\cz_{Z_G(s)^0}\sub\cz_{Z_L(s)^0}$.
\nl
We prove (a).
By 1.4(b) (for $G$ or for $N_GL$), $Z_G(s)$ and $Z_{N_GL}(s)$ are reductive. Since 
$Z_G(s)^0/Q$ is a closed subvariety of the projective variety $G/P$, it 
is itself a projective variety. Hence $Q$ is a parabolic of $Z_G(s)^0$. 
Now $U_P\cap Z_G(s)=U_P\cap Z_G(s)^0$ (see 1.11) is a normal unipotent subgroup of $Q$
hence it is contained in $U_Q$.
Since $L$ has finite index in $N_GL$, $Z_L(s)$ has finite index in $Z_{N_GL}(s)$ hence 
$Z_L(s)$ is reductive. Let $y\in Q$. Then $y\in P$ hence we can write uniquely $y=xu$ 
where $x\in L,u\in U_P$. Since $sys\i=y$ we have $y=(sxs\i)(sus\i)$ and $sxs\i\in L$,
$sus\i\in U_P$. By uniqueness we have $sxs\i=x,sus\i=u$. Thus, $x\in Z_L(s)$,
$u\in U_P\cap Z_G(s)$. The map $y\m x$ is a morphism of algebraic groups 
$f:Q@>>>Z_L(s)$. Since $Q$ is connected, $f(Q)$ is contained in $Z_L(s)^0$. We see that
$Q=Z_L(s)^0(U_P\cap Z_G(s))$. Clearly, $Z_L(s)^0\cap (U_P\cap Z_G(s))=\{1\}$. Since 
$Z_L(s)^0$ is reductive and $U_P\cap Z_G(s)\sub U_Q$, it follows that (a) holds.

We prove (b). Let
$T_1$ be a maximal torus of $Z_L(s)^0$. Let $T_2$ be a maximal torus of $Z_G(s)^0$ 
containg $T_1$. By 1.4(d) (for $N_GL$ instead of $G$), $\tT:=Z_L(T_1)$ is a maximal
torus of $L$. By 1.4(d) (for $G$), $T:=Z_{G^0}(T_2)$ is a maximal torus of $G^0$. Now 
$(\cz_L\cap Z_L(s))^0$ is contained in $\cz_{Z_L(s)^0}$ hence is contained in $T_1$. 
Hence $Z_{G^0}(T_1)\sub Z_{G^0}((\cz_L\cap Z_L(s))^0)=L$ where the last equality comes 
from 1.10(a). Since $T\sub Z_{G^0}(T_1)$, it follows that $T\sub L$ and since 
$T_1\sub T$ and $T$ is commutative, we have $T\sub Z_L(T_1)=\tT$. Since $\tT,T$ are 
maximal tori of $G^0$, we must have $\tT=T$. By 1.4(d) we have 
$T_1=Z_{\tT}(s)^0$, $T_2=Z_T(s)^0$. It follows that $T_1=T_2$. This proves the first 
assertion of (b). Since $Z_G(s)^0$ is reductive, its centre is contained in any maximal
torus of $Z_G(s)^0$, in particular in $T_2$. Since $T_1=T_2$ and $T_1\sub Z_L(s)^0$, we
see that $\cz_{Z_G(s)^0}\sub Z_L(s)^0$ and (b) follows.

\subhead 1.13\endsubhead
Let $P$ be a parabolic of $G^0$. Let $s\in N_GP$ be semisimple. Assume that
$Z_G(s)^0\sub P$. We show:

(a) {\it there is a unique Levi $L$ of $P$ such that $Z_G(s)^0\sub L$. We have 
$s\in N_GL$.}
\nl
Let $T_0$ be a maximal torus of $Z_G(s)^0$. Let $T$ be a maximal torus of $P$ 
containing $T_0$. By 1.4(d), $Z_{G^0}(T_0)$ is a maximal torus of $G^0$. It contains 
$T$ hence it is equal to $T$. Clearly, $Z_{G^0}(T_0)$ is normalized by $s$; hence 
$sTs\i=T$. Let $L$ be the unique Levi of $P$ such that $T\sub L$. Now $sLs\i$ is a Levi
of $sPs\i=P$ containing $sTs\i=T$. By uniqueness, we have $sLs\i=L$. By 1.12(a), 
$U_P\cap Z_G(s)^0$ equals the unipotent radical of $Z_P(s)^0=Z_G(s)^0$ hence it
is $1$. This implies, by 1.12(a), that $Z_P(s)^0$ is a parabolic of $Z_G(s)^0$ with 
Levi $Z_L(s)^0$ and with unipotent radical $\{1\}$. Hence $Z_P(s)^0=Z_G(s)^0=Z_L(s)^0$.
In particular, $Z_G(s)^0\sub L$. This proves the existence of $L$. Assume now that $L'$
is another Levi of $P$ such that $Z_G(s)^0\sub L'$. Since $T_0\sub Z_G(s)^0$, we have 
$T_0\sub L'$. Hence $T_0$ is contained in a maximal torus $T'$ of $L'$. Now $T'$ is 
also a maximal torus of $G^0$. Since $T_0$ is contained in a unique maximal torus of 
$G^0$ (see 1.4(d)) we have $T=T'$. Thus, $L,L'$ are Levi subgroups of $P$ containing a 
common maximal torus of $P$. It follows that $L=L'$. This proves (c).

\subhead 1.14\endsubhead
Let $g\in G$ be quasi-semisimple and let $T_1$ be a maximal torus of $Z_G(g)^0$. Let 
$\cn=\{n\in G^0;ngT_1n\i=gT_1\}$. We show:

(a) $\cn^0=T_1$ {\it hence $\cn/T_1$ is finite; }

(b) {\it any element of $gT_1$ is quasi-semisimple;}

(c) {\it any quasi-semisimple element $g'$ in $gG^0$ is $G^0$-conjugate to some element
in $gT_1$;}

(d) {\it two elements $g',g''\in gT_1$ are in the same $G^0$-conjugacy class if and 
only if they are in the same $\cn/T_1$-orbit on $gT_1$ for the $\cn/T_1$-action induced
by the conjugation action of $\cn$ on $gT_1$. }
\nl
(Closely related results for non connected compact Lie groups appear in \cite{\DS}.)

We prove (a). Let $T$ be the unique maximal torus of $G^0$ such that $T_1\sub T$. We
have $T_1=(T\cap Z_G(g))^0$. (See 1.4(d).) If $n\in\cn$ then $ngn\i=g\t$ with 
$\t\in T_1$ and $gT_1=ngT_1n\i=g\t nT_1n\i$ hence $nT_1n\i=T_1$. Thus, 
$\cn\sub N_{G^0}(T_1)$. Since $\cn^0\sub N_{G^0}(T_1)$ we must have 
$\cn^0\sub Z_{G^0}(T_1)=T$ by a standard rigidity argument. If $n\in\cn^0$ then 
$ngn\i=gt_n$ with $t_n\in T_1$; since $n\in Z_{G^0}(T_1)$, we see that $n\m t_n$ is a 
morphism of algebraic groups $f:\cn^0@>>>T_1$. Clearly, we can find $k\ge 1$ such that 
$g^k$ is in $\cz_{Z_G(s)^0}$ hence in $T_1$ (since $Z_G(s)^0$ is reductive). Then 
$g^k=ng^kn\i=(ngn\i)^k=(gt_n)^k=g^kt_n^k$ hence $t_n^k=1$ for all $n\in\cn^0$. Thus, 
$f(\cn^0)$ is contained in a finite subgroup of $T_1$; being connected, it is $\{1\}$. 
Thus $\cn^0\sub T\cap Z_G(g)$ hence $\cn^0\sub(T\cap Z_G(g))^0=T_1$. The inclusion 
$T_1\sub\cn^0$ is obvious and (a) follows.

We prove (b). The conjugation action of $T_1$ on the variety of all Borels of $G^0$ 
that are normalized by $g$ must have a fixed point since this variety is projective. 
Thus there exists a Borel $B$ of $G^0$ such that $T_1\sub B,gBg\i=B$. Now $T_1$ is
contained in some maximal torus of $B$, which is necessarily $T$, by the definition of
$T$. Since $gTg\i$ is a maximal torus of $G$ that contains $T$, we have $gTg\i=T$
(again by the definition of $T$). If $t\in T_1$ then $t$ normalizes both $B$ and $T$
since $t\in T$. Hence $gt$ normalizes $B$ and $T$. Thus, $gt$ is quasi-semisimple.

We prove (c). Let $T'_1$ be a maximal torus in $Z_G(g')^0$. As in the proof of (b)
we can find a Borel $B'$ of $G^0$ and a maximal torus $T'$ of $B'$ such that
$T'_1\sub T'$ and $g'B'g'{}\i=B',g'T'g'{}\i=T'$. We can find $h\in G^0$ such that
$hB'h\i=B,hT'h\i=T$ (with $B,T$ as above). Let $g''=hg'h\i$. Then $g''Bg''{}\i=B,
g''Tg''{}\i=T$. We also have $gBg\i=B,gTg\i=T$. We have $g''=gy$ where $y\in G^0$
satisfies $yBy\i=B,yTy\i=T$. It follows that $y\in T$. Since a power of $g$ is in 
$T_1$,
$\Ad(g):T@>>>T$ has finite order. Using 1.2(a) we can write $y=g\i y_2gy_2\i y_1$ with
$y_2\in T,y_1\in(T\cap Z_G(g))^0=T_1$. Then $gy=y_2gy_1y_2\i$. We see that $gy$ is
$T$-conjugate (hence $G^0$-conjugate) to $gy_1$. Hence $g'$ is $G^0$-conjugate to 
$gy_1\in gT_1$. This proves (c).

We prove (d). Let $g',g''\in gT_1$ be such that $g''=xg'x\i$ where $x\in G^0$. We have 
$g''=gt$ with $t\in T_1$. Clearly, $T_1$ is a maximal torus of $Z_G(g')^0$ and a
maximal torus of $Z_G(g'')^0$. Then $x\i T_1x$ is a maximal torus of 
$Z_G(x\i g''x)^0=Z_G(g')^0$. The maximal tori $T_1,x\i T_1x$ of $Z_G(g')^0$ are 
conjugate in $Z_G(g')^0$ that is, there exists $z\in Z_G(g')^0$ such that 
$zx\i T_1xz\i=T_1$. Let $n=zx\i$. We have $nT_1n\i=T_1$ and 
$n\i g'n=xz\i g'zx\i=xg'x\i=g''$. We have 
$n\i gT_1n=n\i g'T_1n=g''n\i T_1n=g''T_1=gT_1$ so that $n\i\in\cn$. This proves (c).

\subhead 1.15\endsubhead
We shall need the following result.

(a) {\it The number of unipotent $G$-conjugacy classes of $G$ is finite.
The number of unipotent $G^0$-conjugacy classes of $G$ is finite.}
\nl
These two statements are clearly equivalent. In the case where $G=G^0$ is connected,
(a) is proved in \cite{\FU}. The author handled also the general case by a method
similar to that in \cite{\FU}. See \cite{\SP, 4.1}.

\subhead 1.16\endsubhead
Let $\c\in\Hom(\kk^*,G^0)$. For any $k\in\ZZ$ we set
$\fg_k=\{x\in\fg;\Ad(\c(a))x=a^kx\quad\frl a\in\kk^*\}$. Then 
$\sum_{k\ge 0}\fg_k=\Lie P_\c$ for a well defined parabolic $P_\c$ of $G^0$.
We have $\sum_{k>0}\fg_k=\Lie U_{P_\c}$.

\subhead 1.17\endsubhead
Let $Q$ be a parabolic of $G^0$. Let $g\in N_GQ$. We show that:

(a) {\it there exists $u\in U_Q$ and $\c\in\Hom(\kk^*,G^0)$ such that 
$g\c(a)g\i=u\c(a)u\i$ for all $a\in\kk^*$ and} $P_\c=Q$.
\nl
Let $\p:Q@>>>Q/U_Q$ be the obvious map.
As one easily checks, one can find $\c'\in\Hom(\kk^*,G^0)$ such that 
$\c'(\kk^*)\sub\p\i(\cz_{Q/U_Q}^0)$ and $P_{\c'}=Q$. Let $T$ be a maximal torus of 
$\p\i(\cz_{Q/U_Q}^0)$ that contains $\c'(\kk^*)$. We can find $n\ge 1$ such that 
$g^n\in Q$. For $j\in[0,n-1]$, $g^jTg^{-j}$ is a maximal torus of 
$\p\i(\cz_{Q/U_Q}^0)$, a connected solvable group with unipotent radical $U_Q$. Hence
we can find $u_j\in U_Q$ such that $g^jTg^{-j}=u_jTu_j\i$. Define 
$\c_j\in\Hom(\kk^*,T)$ by $\c_j(a)=u_j\i g^j\c'(a)g^{-j}u_j$ and $\c\in\Hom(\kk^*,T)$ 
by $\c(a)=\c_0(a)\c_1(a)\do\c_{n-1}(a)$. Define an automorphism 
$f:\cz_{Q/U_Q}^0@>>>\cz_{Q/U_Q}^0$ by $f(\p(x))=\p(gxg\i)$ for all
$x\in\p\i(\cz_{Q/U_Q}^0)$. For $a\in\kk^*$ we have 
$\p(\c(a))=\p(\c'(a))f(\p(\c'(a)))f^2(\p(\c'(a)))\do f^{n-1}(\p(\c'(a)))$. Since 
$f^n=1$ it follows that $f(\p(\c(a)))=\p(\c(a))$ that is, $\p(g\c(a)g\i)=\p(\c(a))$. By
a standard argument, if $\l,\l'\in\Hom(\kk^*,Q)$ are such that $\p(\l(a))=\p(\l'(a))$ 
for all $a\in\kk^*$, then there exists $u\in U_Q$ such that $\l(a)=u\l'(a)u\i$ for all 
$a$. Thus, there exists $u\in U_Q$ such that $g\c(a)g\i=u\c(a)u\i$ for all $a$. We have
$P_{\c_j}=u_j\i g^jP_{\c'}g^{-j}u_j=u_j\i g^jQg^{-j}u_j=u_j\i Qu_j=Q$. Hence the 
$\kk^*$-action $a\m\Ad(\c_j(a))$ has $>0$ weights on $\Lie U_Q$ and $<0$ weights on 
$\fg/\Lie Q$. Since these actions (for $j=0,1,\do,n-1$) commute with each other, it 
follows that the $\kk^*$-action 
$a\m\Ad(\c(a))=\Ad(\c_0(a))\Ad(\c_1(a))\do\Ad(\c_{n-1}(a))$ has $>0$ weights on 
$\Lie U_Q$ and $<0$ weights on $\fg/\Lie Q$. Hence $P_\c=Q$. This proves (a).

\subhead 1.18\endsubhead
Let $g\in G$. Let $Q$ be a parabolic subgroup of $Z_G(g_s)^0$ such that $g_uQg_u\i=Q$.
We show that 

(a) {\it there exists a parabolic $P$ of $G^0$ such that $P\cap Z_G(g_s)^0=Q$
and $gPg\i=P$.}
\nl
By 1.17, we can find $\c:\kk^*@>>>Z_G(g_s)^0$ such that $P_\c$ (relative to 
$Z_G(g_s)^0$) is $Q$ and there exists $u\in U_Q$ such that $g_u\c(a)g_u\i=u\c(a)u\i$ 
for all $a\in\kk^*$. Since $g_s\c(a)g_s\i=\c(a)$ for all $a$, we have 
$g\c(a)g\i=u\c(a)u\i$
for all $a$. Define $\ti\c:\kk^*@>>>G^0,\ti\c':\kk^*@>>>G^0$ by $\ti\c(a)=\c(a)$,
$\ti\c'(a)=g\ti\c(a)g\i=u\ti\c(a)u\i$ for all $a$. Let $P=P_{\ti\c}$ (relative to $G$).
From the definition we have $\Lie P\cap\Lie Z_G(g_s)^0=\Lie Q$. Hence 
$P\cap Z_G(g_s)^0=Q$. We have $gPg\i=P_{\ti\c'}=uPu\i=P$ since $u\in U_Q\sub Q\sub P$. 
This proves (a).

\subhead 1.19\endsubhead
Let $T,T'$ be tori and let $f\in\Hom(T',T)$ be surjective. Let $\t:T@>>>T,\t':T'@>>>T'$
be automorphisms of finite order with fixed point
sets $T^\t,T'{}^{\t'}$ such that $f\t'=\t f$. We show that

(a) {\it $f$ restricts to a surjective homomorphism} $(T'{}^{\t'})^0@>>>(T^\t)^0$.
\nl
This can be reduced to the analogous statement where $T,T'$ are replaced by their 
groups of co-characters tensored by $\QQ$. In that case we use the fact that an 
automorphism of finite order of a finite dimensional $\QQ$-vector space is semisimple.

\subhead 1.20\endsubhead
Let $\p:G@>>>G_{ss}$ be the canonical map. Let $a\in G$ be semisimple. Let
$\ph:Z_G(a)^0@>>>Z_{G_{ss}}(\p(a))^0$ be the homomorphism induced by $\p$. We show that

(a) $\ph$ {\it is surjective and} $\Ker(\ph)\sub\cz_{Z_G(a)^0}^0$.

(b) $\ph$ {\it induces a surjective homomorphism }
$\cz_{Z_G(a)^0}^0@>>>\cz_{Z_{G_{ss}}(\p(a))^0}^0$.
\nl
Let $I$ be the image of $\ph$. Let $T$ be a maximal torus of $Z_{G_{ss}}(\p(a))^0$. 
Then $T'=\p\i(T)$ is a torus in $G^0$ and $\Ad(a)(T')=T'$. Moreover, since 
$Z_{G^0}T'$ has finite index in $N_GT'$, we see that there exists an integer $n\ge 1$ 
such that $\Ad(a)^n(t')=t'$ for all $t'\in T'$. Let $T''=\{t'\in T';\Ad(a)(t')=t'\}^0$.
The obvious homomorphism $T'@>>>T$ restricts to a homomorphism $T''@>>>T'$ which is
surjective, by 1.19(a). Since $T''\sub Z_G(a)^0$, we see that $T'\sub I$. Thus, $I$ 
contains the union of all maximal tori of $Z_{G_{ss}}(\p(a))^0$, which is dense in 
$Z_{G_{ss}}(\p(a))^0$, since $Z_{G_{ss}}(\p(a))^0$ is reductive, connected. Thus, $I$ 
is dense in $Z_{G_{ss}}(\p(a))^0$. It is clearly closed, hence $I=Z_{G_{ss}}(\p(a))^0$.
This proves the first assertion of (a). The second assertion of (a) is obvious. Now (b)
is a special case of the following general statement. Let $H@>>>H'$ be a surjective 
homomorphism of connected reductive groups whose kernel is contained in the centre of 
$H$. Then the induced homomorphism $\cz_H^0@>>>\cz_{H'}^0$ is surjective. 

\subhead 1.21\endsubhead
If $D$ is a connected component of $G$, we set

(a) ${}^D\cz_{G^0}=\cz_{G^0}\cap Z_G(g)$
\nl
where $g\in D$. (This does not depend on the choice of $g$.) We write ${}^D\cz_{G^0}^0$
instead of $({}^D\cz_{G^0})^0$. Now let $X$ be a subset of $D$ stable under 
$G^0$-conjugacy. We show:

(b) if ${}^D\cz_{G^0}^0X\sub X$, then $\cz_{G^0}^0X\sub X$.
\nl
(The converse is obvious.) Let $z\in\cz_{G^0}^0,g\in X$. We must show that $zg\in X$. 
Clearly, some power of $g$ is in $G^0$ hence some power of
$\Ad(g):\cz_{G^0}^0@>>>\cz_{G^0}^0$ is $1$. Using 1.2(a), we can write $z=txgxg\i$ with
$t\in(\cz_{G^0}^0\cap Z_G(g))^0={}^D\cz_{G^0}^0,x\in\cz_{G^0}^0$. Then
$zg=txgx\i\in txXx\i=tX=X$. This proves (b).

Consider the $\cz_{G^0}^0\T G^0$-action

(c) $(z,x):y\m xzyx\i$ 
\nl
on $G$ or $D$. This restricts to a ${}^D\cz_{G^0}^0\T G^0$-action on $G$. From (b) we 
see that

(d) {\it The action (c) of $\cz_{G^0}^0\T G^0$ on $D$ and its restriction to
${}^D\cz_{G^0}^0\T G^0$ have exactly the same orbits, that is, any orbit for one action
is an orbit for the other action.}

\subhead 1.22\endsubhead
Let $C$ be an orbit of the $\cz_{G^0}^0\T G^0$-action 1.21(c) on $G$. Let $D$ be the 
connected component of $G$ that contains $C$. Let $C'=\{y\in D;y_s\in C_s\}$. We show:

(a) {\it $C_s$ and $C'$ are closed in $G$ and $y\m y_s$ is a morphism $C'@>>>C_s$};

(b) {if $\bC$ is the closure of $C$ in $G$ and $h\in\bC$ then there exists $h'\in C$
such that $h_s=h'_s$ and $h'{}\i h\in Z_G(h_s)^0$.}
\nl
We prove (a). Let $\s\in C$. By 1.21(d), we have 
$C=\{xz\s x\i;x\in G^0,z\in{}^D\cz_{G^0}\}$. For $x,z$ as above we have 
$(xz\s x\i)_s=xz\s_sx\i$ since $z=z_s,z\s=\s z$. Thus, 
$C_s=\{xz\s_sx\i;x\in G^0,z\in{}^D\cz_{G^0}^0\}$, so that $C_s$ is an orbit
for the action 1.21(c) of ${}^D\cz_{G^0}^0\T G^0$ on $G$. 
We may assume that $G$ is generated by $D$. Then ${}^D\cz_{G^0}^0$ is a closed normal 
subgroup of $G$ and we may form $G'=G/{}^D\cz_{G^0}^0$. Let $\p:G@>>>G'$ be the 
obvious homomorphism. Then $\p(C_s)$ is a semisimple $G'{}^0$-conjugacy class in $G'$
hence $\p(C_s)$ is closed in $G'$ by 1.4(e). Let $G'\sub GL_n(\kk)$ be an imbedding of 
algebraic groups with $n\ge 1$. Let $Y$ be the (semisimple) class in $GL_n(\kk)$ that 
contains $\p(C_s)$. Let $Y'=\{h\in GL_n(\kk);h_s\in Y\}$. It is well known that $Y'$ is
closed in $GL_n(\kk)$ and $Y'@>>>Y,g\m g_s$ is a morphism of varieties. Hence 
$Y'\cap G'$ is closed in $G'$ and $\r:Y'\cap G'@>>>Y\cap G',g\m g_s$ is a morphism of 
varieties. Let $\cx=\{g\in\p(D);g_s\in\p(C_s)\}$. Since $\p(C_s)$ is closed in 
$Y\cap G'$ and 
$\cx=\r\i(\p(C_s))\cap\p(D)$, we see that $\cx$ is closed in $Y'\cap G'$. Next we note
that $C'=\p\i(\cx)$. (The inclusion $C'\sub\p\i(\cx)$ is obvious; the reverse inclusion
follows from the equality ${}^D\cz_{G^0}^0C_s=C_s$.) We see that $C'$ is closed in $D$.
Let $a\in C_s$. Let $H$ be the isotropy group of $a$ in ${}^D\cz_{G^0}^0\T G^0$ (which 
acts transitively on $C_s$). Let $R=\{g\in D;g_s=a\}$. The 
${}^D\cz_{G^0}^0\T G^0$-action on $D$ restricts to an action on $C'$ and this restricts
to an $H$-action on $R$ which induces an isomorphism of algebraic varieties
$({}^D\cz_{G^0}^0\T G^0)\T_HR@>\si>>C'$. Via this isomorphism and the isomorphism 
$({}^D\cz_{G^0}^0\T G^0)\T_H\{a\}@>\si>>C_s$, the map $C'@>>>C_s,g\m g_s$ becomes the 
morphism $({}^D\cz_{G^0}^0\T G^0)\T_HR@>>>({}^D\cz_{G^0}^0\T G^0)\T_H\{a\}$ induced by 
the obvious map $R@>>>\{a\}$. It follows that $C'@>>>C_s,g\m g_s$ is a morphism of 
algebraic varieties. This proves (a).

We prove (b). We have $\bC\sub D$. From (a) we see that $g\m g_s$ is a morphism 
$\bC@>>>C_s$. This morphism is equivariant with respect to the actions 1.21(c) of 
${}^D\cz_{G^0}^0\T G^0$ (transitive on $C_s$). Since $C$ is a dense subset of $\bC$ 
invariant under this action, it follows that for any $a\in C_s$, the intersection of 
$C$ with $\bC_a=\{x\in\bC;x_s=a\}$ is dense in $\bC_a$. Take $a=h_s$. Then any 
connected component $c$ of $\bC_a$ contains some point of $C\cap\bC_a$. Let $c$ be the 
connected component containing $h$. Let $h'\in C\cap c$. Then $h'_s=a$. Since 
$\bC_a\sub Z_G(a)$ we have $c\sub Z_G(a)$. More precisely, $c$ is contained in a 
connected component of $Z_G(a)$. Thus $h'{}\i h\in Z_G(a)^0$. This proves (b).

\subhead 1.23\endsubhead
Let $\boc$ be a $G^0$-conjugacy class in $G$. Let $D$ be the connected component of $G$
that contains $\boc$. Let $\G=\{z\in{}^D\cz_{G^0}^0;z\boc=\boc\}$. We show:

(a) $\G$ {\it is finite.}
\nl
Let $c\in\boc$. If $z\in\G$, then $zc=cz=gcg\i$ hence $zc_s=c_sz=gc_sg\i$ for some 
$g\in G^0$. Let $T_1$ be a maximal torus of $Z_G(c_s)^0$. We have 
${}^D\cz_{G^0}^0\sub T_1$. By 1.4(a), 1.4(d), the $G^0$-conjugacy class of $c_s$ 
intersects $c_sT_1$ in a finite set that is, $F:=\{hc_sh\i;h\in G^0\}\cap c_sT_1$ is 
finite. Since $c_s\G\sub F$, we see that $c_s\G$ is finite and (a) follows.

From (a) we deduce that, for any subgroup $\t$ of ${}^D\cz_{G^0}^0$ we have

(b) $\dim(\t\boc)=\dim\t+\dim\boc$.

\subhead 1.24\endsubhead
Let $\boc,\boc'$ be two $G^0$-conjugacy classes in $G$ contained in the same connected 
component $D$ of $G$. Let $\d,\d'$ be two subtori in $\z:={}^D\cz_{G^0}^0$. We show:

(a) $\d\boc\cap\d'\boc'$ {\it is a finite union of subsets of the form 
$(\d\cap\d')\boc''$ where $\boc''$ are $G^0$-conjugacy classes in $D$.}
\nl
We may assume that $\d\boc\cap\d'\boc\ne\em$. Let $x\in\d\boc\cap\d'\boc'$. Let 
$\boc_1$ be the $G^0$-conjugacy class of $x$. Then 
$\d\boc=\d\boc_1,\d'\boc'=\d'\boc_1$. Thus we may assume that 
$\boc=\boc'$. Let $\G$ be as in 1.23. Then $\ph:\z\T\boc@>>>\z\boc,(a,c)\m ac$ is a 
principal covering with group $\G$ (which acts by $z:(a,c)\m(az,c)$). Since
$\ph(\d\T\boc)=\d\boc,\ph(\d'\T\boc)=\d'\boc$, we have 
$\ph\i(\d\boc)=\cup_{z\in\G}(\d z\T\boc)$, $\ph\i(\d'\boc)=\cup_{z\in\G}\g(\d'z\T\boc)$
and
$$\align&\ph\i(\d\boc\cap\d'\boc)
=\cup_{z\in\G}(\d z\T\boc)\cap\cup_{z'\in\G}(\d'z'\T\boc)\\&=
\cup_{z,z'\in\G}(\d z\T\boc)\cap(\d'z'\T\boc))=\cup_{z,z'\in\G}((\d z\cap\d'z')\T\boc).
\endalign$$
Hence 
$$\d\boc\cap\d'\boc=\ph(\cup_{z,z'\in\G}((\d z\cap\d'z')\T\boc))=
\cup_{z,z'\in\G}((\d z\cap\d'z')\boc).$$
Now $\d z\cap\d'z'$ is either empty or of the form $(\d\cap\d')z_1$ with $z_1\in\z$.
Since $\G$ is finite (see 1.23) we see that 
$\d\boc\cap\d'\boc=\cup_{z_1}(\d\cap\d')z_1\boc$ where $z_1$ runs over a finite subset 
of $\z$. (For any such $z_1$, $z_1\boc$ is a $G^0$-orbit in $D$.) 

\subhead 1.25\endsubhead
Let $P'$ (resp. $P''$) be a parabolic of $G^0$ with Levi $L'$ (resp. $L''$). Assume 
that $L',L''$ contain a common maximal torus. Then $P'\cap P''$ is a connected group 
with Levi $L'\cap L''$. We show:

(a) {\it any $g\in N_GP'\cap N_GP''$ can be written uniquely in the form $g=z\o$ where
$z\in N_GL'\cap N_GP'\cap N_GL''\cap N_GP'',\o\in U_{P'\cap P''}$. Thus, 
$N_GP'\cap N_GP''$ is a semidirect product
$(N_GL'\cap N_GP'\cap N_GL''\cap N_GP'')U_{P'\cap P''}$.}
\nl
Since $g\i$ normalizes $P'\cap P''$, we see that $g\i(L'\cap L'')g$ is a Levi of 
$P'\cap P''$ hence it equals $\o\i(L'\cap L'')\o$ for some $\o\in U_{P\cap P'}$. Then 
$z:=g\o\i$ normalizes $L'\cap L''$ and also $P'$ and $P''$. Now $L'$ is the unique
Levi of $P'$ that contains $L'\cap L''$. Clearly, $zL'z\i$ is again a Levi of $P'$ that
contains $L'\cap L''$. Hence $zL'z\i=L'$ so that $z\in N_GL'\cap N_GP'$. Similarly, 
$z\in N_GL''\cap N_GP''$. Thus, $z\in N_GL'\cap N_GP'\cap N_GL''\cap N_GP''$ and 
$g=z\o$, as required. To prove uniqueness, it is enough to show that
$(N_GL'\cap N_GP'\cap N_GL''\cap N_GP'')\cap U_{P'\cap P''}=\{1\}$.
If $u\in U_{P'\cap P''}$ normalizes $L'$ and $L''$ then, using $u\in P'$, we have
$u\in L'$ and similarly $u\in L''$ hence $u\in L'\cap L''$; but $L'\cap L''$, being a 
Levi of $P'\cap P''$, has trivial intersection with $U_{P'\cap P''}$ hence $u=1$.

We have 

(b) {\it $U_{P'\cap P''}=(L'\cap U_{P''})(U_{P'}\cap P'')$ (semidirect product) and 
also $U_{P'\cap P''}=(L''\cap U_{P'})(U_{P''}\cap P')$ (semidirect product).}

\subhead 1.26\endsubhead
Let $P$ be a parabolic of $G^0$ with Levi $L$. Then 

(a) {\it any $g\in N_GP$ can be written uniquely in the form $g=z\o$ where
$z\in N_GL\cap N_GP,\o\in U_P$. Thus, $N_GP$ is a semidirect product
$(N_GL\cap N_GP)U_P$.}
\nl
This is a special case of 1.25(a) with $P'=P''=P,L'=L''=L$.

\subhead 1.27\endsubhead
Let $H'$ be an algebraic group, let $H$ be a closed subgroup of $H'$ and let $\boc$ be 
a semisimple $H'{}^0$-conjugacy class in $H'$. Then

(a) {\it $H\cap\boc$ is a finite union of (semisimple) $H^0$-conjugacy classes in $H$.}
\nl
We can regard $H'$ as a closed subgroup of $GL_n(\kk)$ for some $n\ge 1$. Let $\boc_1$ 
be the conjugacy class in $GL_n(\kk)$ that contains $\boc$. It is enough to prove (a) 
with $\boc$ replaced by $\boc_1$. Thus we may assume that $H'=GL_n(\kk)$. Let $D$ be a 
connected component of $H$ that contains some semisimple elements. We can find a closed
diagonalizable subgroup $T_D$ of $H$ such that any $H^0$-conjugacy class in 
$D$ meets $T_D$. (We pick a semisimple element $s\in D$ and a maximal torus $T_1$ in 
$Z_H(s)^0$. We can take $T_D$ to be the subgroup generated by $T_1$ and $s$. The fact
that this has the required properties can be deduced from the analogous property of the
reductive group $H/U_{H^0}$, see 1.14(c).) It is enough to prove that $T_D\cap\boc$ is 
finite. Now $T_D$ is contained in a maximal torus $T$ of $GL_n(\kk)$ and it is enough 
to show that $T\cap\boc$ is finite. But this is a well known property of $GL_n(\kk)$.

\head 2. Isolated elements of $G$\endhead
\subhead 2.1\endsubhead
Let $g\in G$. Let
$$\align&T(g)=(\cz_{Z_G(g_s)^0}\cap Z_G(g_u))^0=(\cz_{Z_G(g_s)^0}^0\cap Z_G(g_u))^0
\\&=(\cz_{Z_{G^0}(g_s)^0}\cap Z_{G^0}(g_u))^0
=(\cz_{Z_G(g_s)^0}\cap Z_G(g))^0\endalign$$
(a torus, since $Z_G(g_s)^0$ is reductive). We have $T)g)\sub Z_G(g)$. We sometimes 
write $T_G(g)$ instead of $T(g)$. Clearly,
$$T(xgx\i)=xT(g)x\i \text{ for any } x\in G.$$
Let
$$L(g)=Z_{G^0}(T(g)),\hL(g)=N_G(L(g)).$$
Then 

(a) $L(g)$ {\it is the Levi of a parabolic $P$ of $G^0$ such that $gPg\i=P$.}
\nl
Indeed, we can find $\c\in\Hom(\kk^*,G^0)$ such that $\c(\kk^*)\sub T(g)$ and
$L(g)=Z_{G^0}(\c(\kk^*))$. Then $P=P_\c$, see 1.16, is as required.

From the definition we have

(b) $Z_G(g_s)^0\sub L(g)$. 
\nl
The next result shows that $L(g)$ is characterized by being minimal with the 
properties (a),(b).

(c) {\it Let $Q$ be a parabolic of $G^0$ with Levi $L$ such that
$g\in N_GL\cap N_GQ$ and $Z_G(g_s)^0\sub L$. Then $L(g)\sub L$.}
\nl
Since $Z_{G^0}((\cz_L\cap Z_L(g))^0)=L$, see 1.10(a), it is enough to show that
$Z_{G^0}(T(g))\sub Z_{G^0}((\cz_L\cap Z_L(g))^0)$ or that
$(\cz_L\cap Z_L(g))^0\sub T(g)$ or that $(\cz_L^0\cap Z_G(g))^0\sub T(g)$. Clearly,
$(\cz_L^0\cap Z_G(g))^0\sub Z_G(g_s)$ hence $(\cz_L^0\cap Z_G(g))^0\sub Z_G(g_s)^0$. 
Since $Z_G(g_s)^0\sub L$, we have $\cz_L\cap Z_G(g_s)^0\sub\cz_{Z_G(g_s)^0}$ hence
$(\cz_L^0\cap Z_G(g))^0\sub\cz_{Z_G(g_s)^0}$. Clearly, 
$(\cz_L^0\cap Z_G(g))^0\sub Z_G(g)$ hence
$(\cz_L^0\cap Z_G(g))^0\sub\cz_{Z_G(g_s)^0}\cap Z_G(g)$ and 
$(\cz_L^0\cap Z_G(g))^0\sub(\cz_{Z_G(g_s)^0}\cap Z_G(g))^0=T(g)$, as required.

Since $T(g)\sub Z_G(g)$, we have $g\in\hL(g)$ hence $T_{\hL(g)}(g)$ is defined. We have

(d) $T_{\hL(g)}(g)=T(g)$.
\nl
Since $Z_{\hL(g)}(g_u)=Z_G(g_u)\cap\hL(g)$ and $\cz_{Z_{\hL(g)}(g_s)^0}\sub\hL(g)$, we 
have 
$T_{\hL(g)}(g)=(\cz_{Z_{\hL(g)}(g_s)^0}\cap Z_{\hL(g)}(g_u))^0=
(\cz_{Z_{\hL(g)}(g_s)^0}\cap Z_G(g_u))^0$. From (b) we have 
$Z_G(g_s)^0\sub L(g)\sub\hL(g)$ hence $Z_{\hL(g)}(g_s)^0=Z_G(g_s)^0$ and 
$T_{\hL(g)}(g)=(\cz_{Z_G(g_s)^0}\cap Z_G(g_u))^0=T(g)$. This proves (d). 

We shall need the following result.

(e) {\it  Let $g,g'\in G$ be such that $g_s=g'_s$ and $g'{}\i g\in Z_G(g_s)^0$. Then 
$T(g)=T(g')$.}
\nl
We must show that 
$(\cz_{Z_G(g_s)^0}^0\cap Z_G(g))^0=(\cz_{Z_G(g_s)^0}^0\cap Z_G(g'))^0$. It is enough to
show that, for $x\in\cz_{Z_G(g_s)^0}^0$, the conditions $xg=gx$ and $xg'=g'x$ are 
equivalent. Since $x$ commutes with any element of $Z_G(g_s)^0$, it commutes with
$g'{}\i g$. This proves (e).

\subhead 2.2\endsubhead
We show that the following five conditions for $g\in G$ are equivalent:

(i) $L(g)=G^0$;

(ii) $T(g)\sub\cz_{G^0}$;

(iii) $T(g)={}^D\cz_{G^0}^0$ {\it where $D$ is the connected component of $G$
containing $g$;}

(iv) {\it there is no proper parabolic $P$ of $G^0$ with Levi $L$ such that
$g\in N_GL\cap N_GP,Z_G(g_s)^0\sub L$;}

(v) {\it there is no proper parabolic $P$ of $G^0$ such that
$g\in N_GP,Z_G(g_s)^0\sub P$.}
\nl
Indeed, it is clear that (iii)$\imp$(ii)$\lra$(i). Assume now that (ii) holds. Then any
element of $T(g)$ commutes with any element of $G^0$; since $T(g)\sub Z_G(g)$, we have
$T(g)\sub\cz_{G^0}\cap Z_G(g)$. Since $T(g)$ is connected, we have
$T(g)\sub(\cz_{G^0}\cap Z_G(g))^0$. The reverse inclusion is obvious. We see that (iii)
holds. Thus, the equivalence of (i),(ii),(iii) is established. The equivalence of (i)
and (iv) follows from 2.1(c). It remains to prove the equivalence of (iv),(v). It is
enough to prove the following statement. 

If $P$ is a parabolic of $G^0$ such that $g\in N_GP,Z_G(g_s)^0\sub P$ then there exists
a Levi of $P$ that contains $Z_G(g_s)^0$ and is normalized by $g$.
\nl
By 1.13(a) there is a unique Levi $L$ of $P$ such that $Z_G(g_s)^0\sub L$. Now $gLg\i$ 
is a Levi of $gPg\i=P$ containing $gZ_G(g_s)^0g\i=Z_G(g_s)^0$. By uniqueness, we have 
$gLg\i=L$. This completes the prove of equivalence of (i)-(v).

We say that $g$ is {\it isolated} (in $G$) if it satisfies the equivalent 
conditions (i)-(v). 

If $g\in G$ and $G'$ is a closed subgroup of $G$ containing $g$ and $G^0$ then, 
clearly,

(a) $g$ is isolated in $G$ if and only if $g$ is isolated in $G'$.
\nl
By 2.1(d) we have for $g\in G$:

$Z_{\hL(g)}(T_{\hL(g)}(g))^0=Z_{\hL(g)}(T(g))^0=(\hL(g)\cap Z_G(T(g))^0(\hL(g))^0$,
\nl
hence

(b) {\it $g$ is isolated in $\hL(g)$.}

\subhead 2.3\endsubhead
Let $\p:G@>>>G_{ss}$ be the obvious map. Let $g\in G$. We show:

(a) {\it $g$ is isolated in $G$ if and only if $\p(g)$ is isolated in $G_{ss}$.
Equivalently, the set of isolated elements in $G$ is the inverse image under $\p$ of 
the set of isolated elements in $G_{ss}$.}
\nl
Using the criterion 2.2(v) we see that it is enough to show that conditions (i),(ii) 
below are equivalent:

(i) there exists a proper parabolic $P$ of $G$ such that $Z_G(g_s)^0\sub P,gPg\i=P$

(ii) there exists a proper parabolic $\bP$ of $G_{ss}$ such that 
$Z_{G_{ss}}(\p(g)_s)^0\sub\bP$, $\p(g)\bP\p(g)\i=\bP$.
\nl
Assume that (ii) holds. Let $\bP$ be as in (ii). Then $P=\p\i(\bP)$ is as in (i). 
Conversely, assume that (i) holds. Let $P$ be as in (i). Then $\bP=\p(P)$ is as in (ii)
since $\p$ induces a surjection $Z_G(g_s)^0@>>>Z_{G_{ss}}(\p(g)_s)^0$ (see 1.20(a)).

From (a) we see that 

(b) {\it the set of isolated elements in $G$ is stable under left or right translation 
by} $\cz_{G^0}^0$.

\subhead 2.4\endsubhead
Let $g\in G$. By 1.9(a) (applied to $Z_G(g_s)$ instead of $G$) the set of unipotent 
elements in $Z_G(g_s)^0g_u$ that are quasi-semisimple in $Z_G(g_s)$ is a single 
$Z_G(g_s)^0$-conjugacy class. Let $u$ be an element of this set. According to 1.4(c), 
$g_su$ is quasi-semisimple in $G$. Thus there exist a Borel of $G^0$ and a maximal 
torus of it, both normalized by $g_su$; these are automatically normalized by $u$, 
hence $u$ is quasi-semisimple in $G$. Let $H=Z_G(u)$, a reductive group, by 1.4(b). We
have $g_s\in H$. We show that

(a) {\it $\cz_{Z_H(g_s)^0}$ and its subgroup $\cz_{Z_G(g_s)^0}\cap H$ have 
the same identity component.}
\nl
Let $\tG=Z_G(g_s)$. Let $\tH=Z_{\tG}(u)$ (a reductive group by 1.4(b)). We must prove 
that $\cz_{\tH^0}$ and its subgroup $\cz(\tG^0)\cap\tH$ have the same identity 
component. This follows from 1.7(b) applied to $\tG,u$ instead of $G,u$. (Note that 
1.7(b) is applicable by 1.8(a).)

Next we show:

(b) $(\cz_{G^0}\cap Z_G(g))^0=(\cz_{H^0}\cap Z_G(g))^0$.
\nl
By 1.7(b) (with $g$ replaced by $u$), we have $\cz_{H^0}^0=(\cz_{G^0}\cap H)^0$. (Note 
that 1.7(b) is applicable by 1.8(a).) Since $Z_G(g)\sub H$, the left hand side of (b) 
is 
$$(\cz_{G^0}\cap H\cap Z_G(g))^0=
((\cz_{G^0}\cap H)^0\cap Z_G(g))^0=(\cz_{H^0}^0\cap Z_G(g))^0$$
and this equals the right hand side of (b). This proves (b).

\proclaim{Lemma 2.5} In the setup of 2.4, $g$ is isolated in $G$ if and only if $g_s$ 
is isolated in $H$.
\endproclaim
We have 
$$T_H(g_s)=(\cz_{Z_H(g_s)^0}\cap Z_H(g_s))^0=\cz_{Z_H(g_s)}^0,$$
$$\cz_{Z_G(g_s)^0}\cap H=\cz_{Z_G(g_s)^0}\cap Z_G(g_u)=
\cz_{Z_G(g_s)^0}\cap Z_G(g_s)\cap Z_G(g_u)=\cz_{Z_G(g_s)^0}\cap Z_G(g),$$
since $u\in g_uZ_G(g_s)^0$. Hence $T(g)=(\cz_{Z_G(g_s)^0}\cap H)^0$. Using 2.4(a) we 
deduce that $T(g)=T_H(g_s)$. The condition that $g$ is isolated in $G$ is that 
$T(g)=(\cz_{G^0}\cap Z_G(g))^0$. The condition that $g_s$ is isolated in $H$ is that 
$T_H(g_s)=(\cz_{H^0}\cap Z_H(g_s))^0$ that is, $T_H(g_s)=(\cz_{G^0}\cap Z_G(g))^0$ (see
2.4(b)). These two conditions are equivalent since $T(g)=T_H(g_s)$. The lemma is 
proved.

\proclaim{Lemma 2.6} Assume that $\cz_{G^0}^0=\{1\}$. Then the set of semisimple
elements that are isolated in $G$ is a union of finitely many $G^0$-conjugacy
classes.
\endproclaim
We fix a connected component $D$ of $G$ which contains some semisimple element. Let $s$
be a semisimple element in $D$. Let $T_1$ be a maximal torus of $Z_G(s)^0$. By 1.14(c),
any semisimple element in $D$ is $G^0$-conjugate to some element in $sT_1$. Hence it is
enough to show $sT_1$ contains only finitely many elements that are isolated in $G$. 
There exist finitely many closed connected subgroups $H_1,H_2,\do,H_n$
of $G^0$ such that for any $s'\in sT_1$, $Z_G(s')^0$ is one of $H_1,H_2,\do,H_n$. 
(Indeed, with the notation of 1.5, $\Lie Z_G(s')^0$ is spanned by $\Lie T_1$ and by 
some of the lines $\fg_\a,\a\in R$.) We have $sT_1=\sqc_{i\in[1,n]}X_i$ where 
$X_i=\{s'\in sT_1;Z_G^0(s')=H_i\}$. If $s'\in X_i$, the condition that $s'$ is isolated
in $G$ is that $\cz_{H_i}$ is finite. Thus, either all elements in $X_i$ are isolated 
in $G$ or none is isolated. We may assume that $X_i$ is non-empty and consists of 
isolated
elements if and only if $i\in[1,k]$. Here $k\le n$. We fix $s'_i\in X_i$ for each 
$i\in[1,k]$. If $s'\in X_i$ then $s'\in Z_G(H_i)$ and $s'_i\in Z_G(H_i)$ hence 
$s'_i{}\i s'\in Z_G(H_i)$. Now $s'_i{}\i s'\in T_1\sub H_i$. Hence 
$s'_i{}\i s'\in\cz_{H_i}$. Thus, the set of elements of $sT_1$ that are isolated in $G$
is contained in the finite set $\cup_{i\in[1,k]}s'_i\cz_{H_i}$. The lemma is proved.

\proclaim{Lemma 2.7} The action 1.21(c) of $\cz_{G^0}^0\T G^0$ on $G$ leaves stable
the set of isolated elements in $G$ and has only finitely many orbits there.
\endproclaim
From 2.3(a) we see that the first assertion of the lemma holds and that, to prove the 
second assertion, it is enough to show that the conjugation action of $G_{ss}^0$ on the
set of isolated elements in $G_{ss}$  has only finitely many orbits there. Thus we may
assume that $\cz_{G^0}^0=\{1\}$. Let $D$ be a connected component of $G$. Let $Y$ be 
the set of all elements of $D$ that are isolated in $G$. Let $\tY$ be the set of all 
pairs $(g,u)$ where $g\in D$, $u$ is a unipotent element in $Z_G(g_s)^0g_u$ that is 
quasi-semisimple in $Z_G(g_s)$ and $g_s$ is isolated in $Z_G(u)$. Let $\r:\tY@>>>Y$ be
the first projection. (This is well defined and surjective by Lemma 2.5.) Let 
$(g_0,u_0)\in\tY$. Let $H=Z_G(u_0)$ (a reductive group, by 1.4(b)). Using 1.7(b) (which
is applicable in view of 1.8(a)) we see that $\cz_{H^0}^0=\{1\}$. Applying Lemma 2.6 to
$H$ instead of $G$, we can find isolated semisimple elements $s_1,s_2,\do,s_n$ in $H$ 
such that any isolated semisimple element in $H$ is $H^0$-conjugate to some $s_i$. By 
1.15(a), for any $i\in[1,n]$ we can find unipotent elements 
$v_{i1},v_{i2},\do,v_{ip_i}$ in $Z_G(s_i)$ such that any unipotent element in 
$Z_G(s_i)$ is $Z_G(s_i)^0$-conjugate to some $v_{ij}$. It is enough to show:

(a) {\it Let $\co$ be an orbit for the $G^0$-action on $\tY$ given by conjugation on 
both coordinates. Then $(s_iv_{ij},u')\in\co$ for some $i\in[1,n],j\in[1,p_i]$ and some
$u'$.}
\nl
(Indeed, since $\r$ is surjective, (a) would imply that any $G^0$-conjugacy class in
$Y$ contains $s_iv_{ij}$ for some $i\in[1,n],j\in[1,p_i]$.) Let $(g,u)\in\co$. 
Now $uG^0=g_uG^0=(g_0)_uG^0=u_0G^0$ (since $gG^0=g_0G^0=D$). Thus $u,u_0$ are 
unipotent, quasi-semisimple in $G$, in the same connected component of $G$, hence 
$u_0=huh\i$ for some $h\in G^0$. (See 1.9(a).) We have $(hgh\i,huh\i)\in\co$. Setting 
$g'=hgh\i$ we have $(g',u_0)\in\co$. Since $g'_s$ is an isolated semisimple element in 
$H$, we can find $h'\in H^0$ such that $h'g'_sh'{}\i=s_i$ with $i\in[1,n]$. We have 
$(h'g'h'{}\i,h'u_0h'{}\i)\in\co$. Now $h'u_0h'{}\i=u_0$, the semisimple part of
$h'g'h'{}\i$ is $h'g'_sh\i=s_i$ and the unipotent part of $h'g'h'{}\i$ is a unipotent 
element $\ti v_i$ in $Z_G(s_i)$. Thus we have $(s_i\ti v_i,u_0)\in\co$. We can find 
$h''\in Z_G(s_i)^0$ such that $h''\ti v_ih''{}\i=v_{ij}$ with $j\in[1,p_i]$. We have 
$(s_iv_{ij},h''u_0h''{}\i)=(h''s_i\ti v_ih''{}\i,h''u_0h''{}\i)\in\co$. This proves 
(a). The lemma is proved.

\proclaim{Lemma 2.8} The set of isolated elements of $G$ is closed in $G$.
\endproclaim
Using 2.3(a) we see that if the lemma holds for $G_{ss}$ then it holds for $G$. Hence
we may assume that $\cz_{G^0}^0=\{1\}$. Let $G'$ be the set of isolated elements of
$G$. Since $G'$ is a union of finitely many $G^0$-conjugacy classes (see Lemma 2.7),
we see that $G'_s$ is a union of finitely many semisimple $G^0$-conjugacy classes 
$E^1,E^2,\do,E^n$. For $i\in[1,n]$ let $G'{}^i$ be the inverse image of $E^i$ under 
$G'@>>>G'_s,g\m g_s$. It is enough to prove that $G'{}^i$ is closed in $G$ for any 
$i\in[1,n]$. 
Let $E'{}^i=\{g\in G;g_s\in E^i\}$. By 1.22(a), $E^i,E'{}^i$ are closed in $G$ and 
$\p:E'{}^i@>>>E^i,g\m g_s$ is a morphism. Note that $\p$ commutes with the conjugation 
action of $G^0$ on $E^i,E'{}^i$ and that action is transitive on $E^i$. Hence a 
$G^0$-stable subset of $E'{}^i$ is closed if and only if its intersection with 
$\p\i(s)$ is closed in $G$ for some/any $s\in E^i$. Thus, to prove that $G'{}^i$ is 
closed in $E'{}^i$ (hence in $G$) it is enough to prove that, if $s\in E^i$, then 
$\{g\in G'{}^i;g_s=s\}=\{g\in G';g_s=s\}$ is closed in $G$. Let $\t=\cz_{Z_G(s)^0}^0$. 
We must show that $\{su;u\in G,u\text{ unipotent},us=su,(\t\cap Z_G(u))^0=\{1\}\}$
is closed in $G$ or that $\{u\in Z_G(s);u\text{ unipotent},(\t\cap Z_G(u))^0=\{1\}\}$
is closed in $Z_G(s)$. Since the set of unipotent elements in $Z_G(s)$ is closed in
$Z_G(s)$, it is enough to show that $\{g\in Z_G(s);(\t\cap Z_G(g))^0=\{1\}\}$ is closed
in $Z_G(s)$. Let $X$ be a connected component of $Z_G(s)$. It is enough to show that 
$X_0:=\{g\in X;(\t\cap Z_G(g))^0=\{1\}\}$ is closed in $X$. Now $\t\cap Z_G(g)$ depends
only on the connected component of $g$ in $Z_G(s)$. Thus, either $X_0=X$ or $X_0=\em$. 
In both cases, $X_0$ is closed in $X$. The lemma is proved.

\head 3. A stratification of $G$\endhead
\subhead 3.1\endsubhead
For $g,g'\in G$ we write $g\si g'$ if $g'g\i\in T(g)=T(g')$. This is an equivalence 
relation on $G$. For $x\in G^0,g,g'\in G$ with $g\si g'$ we have $xgx\i\si xg'x\i$
since $xT(g)x\i=T(xgx\i)$. Hence the relation on $G$ given by $g\app g'$ if 
$xgx\i\si g'$ for some $x\in G^0$ is an equivalence relation. The equivalence classes 
on $G$ for $\app$ are called the {\it strata} of $G$. Each stratum of $G$ is contained 
in a connected component of $G$ and is stable under conjugation by $G^0$.

\proclaim{Lemma 3.2} Let $\p:G@>>>G_{ss}$ be the obvious homomorphism. Let 
$g,g'\in G$. We have $g\app g'$ (in $G$) if and only if $\p(g)\app\p(g')$ (in 
$G_{ss}$). Thus, the strata of $G$ are exactly the inverse images under $\p$ of the 
strata of $G_{ss}$ and each stratum of $G$ is stable under left or right translation by
$\cz_{G^0}^0$.
\endproclaim
Let $g\in G$. By 1.20(b), $\p$ induces a surjective homomorphism of tori
$\cz_{Z_G(g_s)^0}^0@>>>\cz_{Z_{G_{ss}}(\p(g)_s)^0}^0$. This is compatible with the
automorphisms
$$\Ad(g_u):\cz_{Z_G(g_s)^0}^0@>>>\cz_{Z_G(g_s)^0}^0, 
\Ad(\p(g)_u):\cz_{Z_{G_{ss}}(\p(g)_s)^0}^0@>>>\cz_{Z_G(g_s)^0}^0,$$
which have finite order. (Some power of $g_u$ belongs to $Z_G(g_s)^0$.)
Applying 1.19(a) we see that $\p$ restricts to a surjective homomorphism
$$(\cz_{Z_G(g_s)^0}^0\cap Z_G(g_u))^0@>>>
(\cz_{Z_{G_{ss}}(\p(g)_s)^0}^0\cap Z_{G_{ss}}(\p(g)_u))^0.$$
Thus, we have

(a) $\p(T_G(g))=T_{G_{ss}}(\p(g))$.
\nl
We have $T_G(g)\sub((T_G(g)\cz_{G^0}^0)\cap Z_G(g))^0=T_G(g)(\cz_{G^0}^0\cap Z_G(g))^0
=T_G(g)$. Hence 

(b) $T_G(g)=((T_G(g)\cz_{G^0}^0)\cap Z_G(g))^0$.
\nl
If $z''\in(\cz_{G^0}\cap Z_G(g))^0$ we have

(c) $T_G(z''g)=T_G(g)$ and $z''g\si g$.
\nl
The first assertion follows from 2.1(e) since $(z''g)_s=g_s$ and 
$(z''g)\i g\in Z_G(g_s)^0$. The second assertion follows from the first since 
$z''gg\i=z''\in T(g)$. 

If $z\in\cz_{G^0}^0$ then

(d) $T_G(zg)=T_G(g)$ and $zg\app g'$.
\nl
By 1.2(a) we have $z=z'z''gz'{}\i g\i$ for some 
$z'\in\cz_{G^0}^0,z''\in(\cz_{G^0}^0\cap Z_G(g))^0$. Hence
$$T_G(zg)=T_G(z'z''gz'{}\i)=z'T_G(z''g)z'{}\i=z'T_G(g)z'{}\i=T_G(g)$$
where the third equality comes from (c). We also have 
$zg=z'z''gz'{}\i\app z''g\app g$. (See (c).) This proves (d).

Next, for $g,g'\in G$,

(e) {\it we have $\p(g)\si\p(g')$ (in $G_{ss}$) if and only if $zg\si g'$ (in $G$) for some}
$z\in\cz_{G^0}^0$.
\nl
Assume that $z\in\cz_{G^0}^0$ and $zg\si g'$ that is, $g'g\i z\i\in T_G(zg)=T_G(g')$.
Applying $\p$ and using (a), we obtain 
$\p(g')\p(g)\i\in T_{G_{ss}}(\p(g))=T_{G_{ss}}(\p(g'))$. 
Thus, $\p(g)\si\p(g')$ in $G_{ss}$. 

Conversely, assume that $\p(g)\si\p(g')$ in $G_{ss}$ that is,
$\p(g')\p(g)\i\in T_{G_{ss}}(\p(g))=T_{G_{ss}}(\p(g'))$. Using 
(a), we deduce that $g'g\i\in T_G(g)\cz_{G^0}^0=T_G(g')\cz_{G^0}^0$.
Thus $g'=tzg$ where $t\in T_G(g),z\in\cz_{G^0}^0$. Using (b), we have
$$\align&T_G(g')=((T_G(g')\cz_{G^0}^0)\cap Z_G(g'))^0\\&=((T_G(g)\cz_{G^0}^0)
\cap Z_G(tzg))^0=((T_G(g)\cz_{G^0}^0)\cap Z_G(g))^0=T_G(g)=T_G(zg).\endalign$$
(The third equality holds since $tz$ belongs to the commutative group 
$T_G(g)\cz_{G^0}^0$ and the fifth equality comes from (d).) Thus, $T_G(g')=T_G(zg)$. 
Now $g'(zg)\i=t\in T_G(g)=T_G(zg)$. Hence $g'\si zg$. This proves (e).

The lemma follows from (d),(e).

\subhead 3.3\endsubhead
The set of isolated elements of $G$ is a union of strata of $G$. (Assume that $g\in G$
is isolated and $g'\in G$ is in the stratum of $g$. We must show that $g'$ is isolated.
We may assume that $g'\si g$. Then $T(g')=T(g)$. By assumption, $T(g)\sub\cz_{G^0}^0$.
Hence $T'(g)\sub\cz_{G^0}^0$ and $g'$ is isolated.) The strata of $G$ that consist of
isolated elements are called {\it isolated strata}. We show:

(a) {\it any isolated stratum of $G$ is a single orbit for the action 1.21(c) of
$\cz_{G^0}^0\T G^0$ on $G$.}
\nl
In view of 2.3(a) and Lemma 3.2, it is enough to consider the case where 
$\cz_{G^0}^0=\{1\}$. In this case we must show that any isolated stratum $C$ of $G$ is 
a single $G^0$-conjugacy class in $G$. Let $g,g'\in C$. We have $g'=hzgh\i$ where 
$h\in G^0,z\in T(g)$. Since $g$ is isolated, we have $T(g)=\{1\}$. Hence $g'=hgh\i$. 
This proves (a).

\proclaim{Lemma 3.4} Let $C$ be a stratum of $G$. If $g,g'\in C$ then 
$Z_G(g)^0,Z_G(g')^0$ are $G^0$-conjugate.
\endproclaim
We may assume that $g\si g'$. We have $g'=tg$ where $t\in T(g)$. If $x\in Z_G(g)^0$ 
then $x\in Z_G(g_s)^0$ and $t\in\cz_{Z_G(g_s)^0}$ hence $xt=tx$. We have also $xg=gx$ 
hence $xtg=tgx$. Thus, $x\in Z_G(g')$. We see that $Z_G(g)^0\sub Z_G(g')$ so that 
$Z_G(g)^0\sub Z_G(g')^0$. Interchanging the roles of $g,g'$ we obtain 
$Z_G(g')^0\sub Z_G(g)^0$ hence $Z_G(g)^0=Z_G(g')^0$. The lemma is proved.

\subhead 3.5\endsubhead
Let $\AA$ be the set of all pairs $(L,S)$ where $L$ is a Levi of some parabolic of 
$G^0$ and $S$ is an isolated stratum of $N_GL$ with the following property: there 
exists a parabolic $P$ of $G^0$ with Levi $L$ such that $S\sub N_GP$. Now $G^0$ acts on
$\AA$ by conjugation. Let $G^0\bsl\AA$ be the set of orbits of this action.

\proclaim{Lemma 3.6} Let $C$ be a stratum of $G$. For $g\in C$ let $L=L(g)$ and let 
$S$ be the stratum of $N_GL$ that contains $g$. Then $(L,S)\in\AA$ and $C\m(L,S)$ is a 
well defined injective map from the set of strata of $G$ to $G^0\bsl\AA$. 
\endproclaim
By 2.2(b), $g$ is isolated in $N_GL$. By 2.1(a) there exists a
parabolic $P$ of $G^0$ with Levi $L$ such that $g\in N_GP$. We have $S\sub N_GP$. 
(Since $S$ is an isolated stratum of $N_GL$, any element of $S$ is of the form 
$hzgh\i$ with $z\in\cz_L^0,h\in L$ (see 3.3(a)); clearly, $hzgh\i\in N_GP$.) Thus 
$(L,S)\in\AA$. Let $g'\in G$ be such that $g\app g'$. Let $L'=L(g')$ and let $S'$ be 
the stratum of $N_GL'$ that contains $g'$. We show that $(L,S),(L',S')$ are in the same
$G^0$-orbit. Replacing $g$ by a $G^0$-conjugate, we may assume that $g\si g'$. Then 
$T(g)=T(g')$ hence $L=L'$. Also $g'g\i\in T(g)$. Using 2.1(d) we have 
$T_{N_GL}(g)=T(g)=T(g')=T_{N_GL}(g')$ and $g'g\i\in T_{N_GL}$ hence $g\si g'$ (relative
to $N_GL$). Hence $S=S'$. Thus the map in the lemma is well defined. We show that it is
injective. Assume that $g,g'\in G$ are such that for some $h\in G^0$ we have $hLh\i=L'$
where $L=L(g),L'=L(g')$ and $hgh\i\app g'$ where $\app$ is relative to $N_GL'$. We must
show that $g\app g'$ (relative to $G$). Replacing $g$ by $hgh\i$, we may assume that 
$L=L',g\app g'$ (relative to $N_GL$). Replacing $g$ by an $L$-conjugate, we may assume 
further that $g\si g'$ (relative to $N_GL$). We have 
$g'g\i\in T_{N_GL}(g)=T_{N_GL}(g')$. But $T_{N_GL}(g)=T(g)$, $T_{N_GL}(g')=T(g')$, so 
that $g\si g'$ (relative to $G$). The lemma is proved.

\proclaim{Proposition 3.7} The number of strata of $G$ is finite.
\endproclaim
By Lemma 3.6, it is enough to prove that $G^0\bsl\AA$ is finite. Since the Levi 
subgroups $L$ of parabolics of $G^0$ fall into finitely many $G^0$-conjugacy classes, 
it is enough to show that for any such $L$, there are only finitely many isolated 
strata of $N_GL$. This follows from Lemma 2.7 and 3.3(a). The proposition is proved.

\subhead 3.8\endsubhead
Let $P$ be a parabolic of $G^0$ with Levi $L$. Let $g$ be an isolated element of 
$N_GL\cap N_GP$. (Equivalently, $g$ is an element of $N_GL\cap N_GP$ that is isolated
in $N_GL$, see 2.2(a); indeed, $N_GL\cap N_GP\sub N_GL$ have the same identity 
component, $L$.) We show:

(a) $T(g)\sub(\cz_L\cap Z_L(g))^0$;

(b) $L\sub L(g)$;
\nl
We prove (a). Since  $g_s\in N_GL\cap N_GP$ we see from 1.12(b) that
$\cz_{Z_G(g_s)^0}^0\sub\cz_{Z_L(g_s)^0}^0$. Hence 
$$\align&T(g)=(\cz_{Z_G(g_s)^0}^0\cap Z_G(g))^0\sub(\cz_{Z_L(g_s)^0}^0\cap Z_G(g))^0
\\&=(\cz_{Z_L(g_s)^0}^0\cap Z_L(g))^0=(\cz_L\cap Z_L(g))^0\endalign$$
where the last equality holds since $g$ is isolated in $N_GL$. This proves (a).

From (a) we deduce $Z_{G^0}(T(g))\sps Z_{G^0}((\cz_L\cap Z_L(g))^0)=L$ (the last 
equality comes from 1.10(a)). This proves (b).

\subhead 3.9\endsubhead
Let $P,L,g$ be as in 3.8. We show that the following three conditions are equivalent:

(i) $L=L(g)$; 

(ii) $T(g)=(\cz_L\cap Z_L(g))^0$; 

(iii) $Z_G(g_s)^0\sub L$.
\nl
If (ii) holds, then $Z_{G^0}(T(g))=Z_{G^0}((\cz_L\cap Z_L(g))^0)=L$ (see 1.10(a)) so 
that (i) holds. Assume now that (i) holds. We show that (ii) holds. By 3.8(a) it is 
enough to show that $(\cz_L\cap Z_G(g))^0\sub T(g)$ or that 
$(\cz_{Z_{G^0}(T(g))}\cap Z_G(g))^0\sub T(g)$ or that
$(\cz_{Z_{G^0}(T(g))}\cap Z_G(g_s)^0\cap Z_G(g_u))^0\sub\cz_{Z_G(g_s)^0}\cap Z_G(g_u)$.
It is enough to show that $\cz_{Z_{G^0}(T(g))}\cap Z_G(g_s)^0\sub\cz_{Z_G(g_s)^0}$. This
follows from $Z_G(g_s)^0\sub Z_{G^0}(T(g))$ (a consequence of the definition of 
$T(g)$). 

By 3.8(a), condition (ii) is equivalent to $(\cz_L\cap Z_G(g))^0\sub T(g)$ and also to 
$(\cz_L\cap Z_G(g))^0\sub\cz_{Z_G(g_s)^0}$. Since $(\cz_L\cap Z_G(g))^0\sub Z_G(g_s)^0$
this is also equivalent to the condition $Z_G(g_s)^0\sub Z_{G^0}((\cz_L\cap Z_G(g))^0)$
which by 1.10(a) is the same as (iii). 

\subhead 3.10\endsubhead
Let $P,L,g$ be as in 3.8. For any $z\in(\cz_L\cap Z_L(g))^0$, the element $zg$ is 
isolated in $N_GL\cap N_GP$. (This follows from 2.3(a) for $N_GL\cap N_GL$ instead of 
$G$.) We show:

(a) $\{z\in(\cz_L\cap Z_L(g))^0;Z_G(zg_s)^0\sub L\}$ {\it is open dense in} 
$(\cz_L\cap Z_L(g))^0$.
\nl
Assume first that $g=s$ is semisimple. Since $s$ is isolated in $N_GL$ we have 
$(\cz_L\cap Z_L(s))^0=\cz_{Z_L(s)^0}^0$. Hence it is enough to show:

(b) $\{z\in\cz_{Z_L(s)^0}^0;Z_L(zs)^0=Z_G(zs)^0\}$ {\it is open dense in}
$\cz_{Z_L(s)^0}^0$.
\nl
By 1.12(b) we can find a maximal torus $T_1$ of $Z_L(s)^0$ which is also a maximal 
torus of $Z_G(s)^0$. Define $\fg_\a,R$ as in 1.5 in terms of $s,G,T_1$. We have 
$\cz_{Z_L(s)^0}^0\sub T_1$. For $z\in\cz_{Z_L(s)^0}^0$, $Z_G(zs)^0$ contains $T_1$ and 
is normalized by $s$ hence $\Lie Z_G(zs)^0$ is spanned by $\Lie T_1$ and by the 
$\fg_\a$ with $\a\in R,\a(zs)=1$. Similarly, $\Lie Z_L(zs)^0$ is spanned by $\Lie T_1$ 
and by the $\fg_\a$ with $\a\in R,\a(zs)=1,\fg_\a\sub\Lie L$. For each $\a\in R$ we set
$X_\a=\{z\in\cz_{Z_L(s)^0}^0;\a(zs)\ne 1\}$. We see that it is enough to show: 

(c) $\cap_{\a\in R;\fg_\a\not\sub\Lie L}X_\a$ {\it is open dense in}
$\cz_{Z_L(s)^0}^0$.
\nl
Assume that $\fg_\a\not\sub\Lie L$ but $\a(zs)=1$ for all $z\in\cz_{Z_L(s)^0)^0}$. Then
$\a(s)=1$ and $\a|_{\cz_{Z_L(s)^0}^0}=1$. Thus, 
$\fg_\a\sub\Lie Z_{G^0}(\cz_{Z_L(s)^0}^0)=\Lie Z_{G^0}((\cz_L\cap Z_L(s))^0)=\Lie L$ 
(the last equality comes from 1.10(a)). Now $\fg_a\sub\Lie L$ is a contradiction.
Thus, for any $\a\in R$ such that $\fg_\a\not\sub\Lie L$, the open subset $X_\a$ of 
$\cz_{Z_L(s)^0}^0$ is non-empty hence dense. It follows that (c) holds; (b) is proved.

We now consider the general case. Let $u$ be a unipotent element in $Z_L(g_s)^0g_u$ 
which is quasi-semisimple in $Z_L(g_s)$. As in 2.4, $u$ is quasi-semisimple in $N_GL$. 
Hence there exists a Borel $\b$ of $L$ and a maximal torus $T$ of $\b$ such that $\b$ 
and $T$ are normalized by $u$. Now $U_P$ is normalized by $g$ hence by $g_u$, hence by 
$lg_u$ for any $l\in L$. Since $u=lg_u$ for some $l\in L$, we have $u\in N_G(U_P)$. 
Then $u$ normalizes $\b U_P$, a Borel of $G^0$ containing $T$. Thus, $u$ is 
quasi-semisimple in $G$ and $H=Z_G(u)$ is a reductive group. Similarly, $L_1=Z_L(u)$ is
a reductive group; it is a Levi of the parabolic $P_1=Z_P(u)^0$ of $H^0$. By Lemma 2.5
(with $G$ replaced by $N_GL$), $g_s$ is isolated in $N_H(L_1)$. Also, 
$g_s\in N_H(L_1)\cap N_H(P_1)$. Applying the already proved part of (a) to 
$H,L_1,P_1,g_s$ instead of $G,L,P,g$ we see that 

$\{z\in(\cz_{L_1}\cap Z_{L_1}(g_s))^0;Z_H(zg_s)^0\sub L_1\}$ {\it is open dense in} 
$(\cz_{L_1}\cap Z_{L_1}(g_s))^0$.
\nl
It is then enough to show that

(d) $(\cz_{L_1}\cap Z_{L_1}(g_s))^0=(\cz_L\cap Z_L(g))^0$
\nl
and that, for $z$ in (d), we have $Z_H(zg_s)^0\sub L_1$ if and only if 
$Z_G(zg_s)^0\sub L$, or equivalently (see 3.9), that we have 
$T_H(zg_s)=(\cz_{L_1}\cap Z_{L_1}(zg_s))^0$ if and only if 
$T_G(zg)=(\cz_L\cap Z_L(zg))^0$.
It is enough to prove (d) and that for any $z$ in (d) we have

(e) $T_H(zg_s)=T_G(zg)$,

(f) $(\cz_{L_1}\cap Z_H(zg_s))^0=(\cz_L\cap Z_L(zg))^0$.
\nl
As a special case of 1.7(b) we have $\cz_{L_1}^0=(\cz_L\cap Z_L(u))^0$ hence

$(\cz_{L_1}\cap Z_{L_1}(g_s))^0=(\cz_L\cap Z_L(u)\cap Z_L(ug_s))^0=
(\cz_L\cap Z_L(g))^0,$
\nl
so that (d) holds. (In the last equality we have used that $u\in Z_L(g_s)^0g_u$.)
Similarly, for $z$ in (d) we have

$(\cz_{L_1}\cap Z_{L_1}(zg_s))^0=(\cz_L\cap Z_L(u)\cap Z_L(uzg_s))^0=
(\cz_L\cap Z_L(zg))^0,$
\nl
so that (f) holds. Now (e) is shown as in the proof of Lemma 2.5 (for $zg$ instead of 
$g$). The lemma is proved.

\proclaim{Lemma 3.11} Let $(L,S)\in\AA$. Let $P$ be a parabolic of $G^0$ with Levi $L$ 
such that $S\sub N_GP$. Let $S^*=\{g\in S;Z_G(g_s)^0\sub L\}$. Then $S^*$ is an open 
dense subset of $S$.
\endproclaim   
By 3.3(a) (for $L$), $S$ is contained in a connected component $\d$ of $N_GL\cap N_GP$.
We set ${}^\d\cz_L=\cz_L\cap Z_L(g')$ for some/any $g'\in\d$. By 3.3(a) and 1.21(d), the
action 1.21(c) of ${}^\d\cz_L^0\T L$ on $S$ is transitive. The restriction of this
action to ${}^\d\cz_L^0$ is a free action (left translation) and 
$\hat S={}^\d\cz_L^0\bsl S$ is well defined. The conjugation action of $L$ on $S$
induces a transitive $L$-action on $\hat S$. Hence the condition that an $L$-invariant 
subset $X$ of $S$ is open dense in $S$ is equivalent to the condition that for any 
${}^\d\cz_L^0$-orbit $\co$ on $S$, the intersection $X\cap\co$ is open dense in $\co$. 
For $X=S^*$ this last condition holds by 3.10(a). The lemma is proved.

\proclaim{Proposition 3.12} (a) If $(L,S)\in\AA$ then $Y_{L,S}=\cup_{x\in G^0}xS^*x\i$ 
is a stratum of $G$.

(b) $(L,S)\m Y_{L,S}$ is a bijection between $G^0\bsl\AA$ and the set of strata of $G$.
In particular, we have a partition $G=\sqc_{L,S}Y_{L,S}$ where $L,S$ runs through a set
of representatives for the $G^0$-orbits in $\AA$; this is the same as the partition of 
$G$ into strata.
\endproclaim
We prove (a). By Lemma 3.11, $S^*$ is non-empty. Hence $Y_{L,S}$ is non-empty. Now
$Y_{L,S}$ is contained in a stratum of $G$. (It is enough to show that $S^*$ is
contained in a stratum of $G$. It is also enough to show that, if $g,g'\in S^*$ and
$g\si g'$ relative to $N_GL$, then $g\si g'$ relative to $G$. This follows from the 
equalities $T_G(g)=T_{N_GL(g)}(g)=T_{N_GL}(g),T_G(g')=T_{N_GL(g')}(g')=T_{N_GL}(g')$,
see 2.1(d), 3.9.) We show that $Y_{L,S}$ is a stratum of $G$. It is enough to show 
that, if $g\in S^*$ and $g'\in G,g\app g'$ relative to $G$, then $g'\in Y_{L,S}$. We 
may assume that $g'\si g$ relative to $G$. Since $T_G(g)=T_G(g')$, we have 
$L(g)=Z_{G^0}T_G(g)=Z_{G^0}T_G(g')=L(g')$. Since $g\in S^*$ we have $L=L(g)$ (see 
3.9). It follows that $L(g')=L$. In particular, $g'\in N_GL$. As above, we have 
$T_G(g)=T_{N_GL}(g),T_G(g')=T_{N_GL(g')}(g')$ and the last group equals $T_{N_GL}(g')$ 
since $L(g')=L$. Since $g'{}\i g\in T_G(g)=T_G(g')$, we see that
$g'{}\i g\in T_{N_GL}(g)=T_{N_GL}(g')$. Thus $g'\si g$ relative to $N_LG$. It follows
that $g'\in S$. More precisely, since $L=L(g')$, we have $g'\in S^*$, as required. This
proves (a).

From the definitions it is clear that the map in (b) is the inverse of the map 
$C\m(L,S)$ in Lemma 3.6. The lemma is proved.

\subhead 3.13\endsubhead
Let $(L,S)\in\AA$. Let 
$$\tY_{L,S}=\{(g,xL)\in G\T G^0/L;x\i gx\in S^*\}.$$
Define $\p:\tY_{L,S}@>>>Y_{L,S}$ by $\p(g,xL)=g$. Now
$\cw_S=\{n\in N_{G^0}L;nSn\i=S\}/L$ (a subgroup of the finite group $N_{G^0}L/L$) acts
(freely) on $\tY_{L,S}$ by $n:(g,xL)\m(g,xn\i L)$. 

(a) {\it This makes $\p:\tY_{L,S}@>>>Y_{L,S}$ into a principal $\cw_S$-bundle.}
\nl
We must show that, if $g\in G,x\in G^0,x'\in G^0$ satisfy 
$x\i gx\in S^*,x'{}\i gx'\in S^*$, then $x'=xn\i$ with $n\in N_{G^0}L,nSn\i=S$. 
Replacing $g,x,x'$ by $x'{}\i gx',x'{}\i x,1,$ we see that we may assume that $x'=1$ and 
we must show that $xLx\i=L,xSx\i=S$. We have $g\in S^*,x\i gx\in S^*$ and 
$L=L(x\i gx)=x\i L(g)x=x\i Lx$. If $g'\in S$ then $g'=hzgh\i$ for some $h\in L$, 
$z\in(\cz_L\cap Z_L(g))^0$, see 3.3(a) and 1.21(d). We have 
$$x\i zx\in(\cz_{x\i Lx}\cap Z_{x\i Lx}(x\i gx))^0=
(\cz_L\cap Z_L(x\i gx))^0=(\cz_L\cap Z_L(g))^0$$
(since $x\i gx\in gL$) hence $x\i g'x=(x\i hx)(x\i zx)(x\i gx)(x\i h\i x)\in S$. Thus, 
$x\i Sx\sub S$. Since $S$ is irreducible we have $x\i Sx=S$, as required.

Since $\tY_{L,S}$ is an irreducible variety of dimension

$\dim(G^0/L)+\dim S^*=\dim(G^0/L)+\dim S$, 
\nl
it follows that 

(b) {\it $Y_{L,S}$ is an irreducible constructible subset of $G$ of dimension} 
$\dim(G^0/L)+\dim S$.

\proclaim{Lemma 3.14}  Let $(L,S)\in\AA$. Let $P$ be a parabolic of $G^0$ with Levi $L$ 
such that $S\sub N_GP$. Let $\bS$ be the closure of $S$ in $N_GL$. Let
$Y'_{L,S}=\cup_{x\in G^0}x\bS U_Px\i$. Then the closure of $Y_{L,S}$ in $G$ is 
$Y'_{L,S}$. In particular, $Y'_{L,S}$ is independent of the choice of $P$.
\endproclaim
Let $X=\{(g,xP)\in G\T G^0/P;x\i gx\in\bS U_P\}$. Let $\ph:X@>>>G$ be the first
projection. Then $\ph$ is proper since $G^0/P$ is complete and $\bS U_P$ is a closed 
$\Ad(P)$-stable subset of $N_GP$. Hence $\ph(X)=Y'_{L,S}$ is a closed subset of $G$.
Since $X$ is irreducible of dimension $\dim(G^0/P)+\dim S+\dim U_P=\dim(G^0/L)+\dim S$,
we see that $Y'_{L,S}$ is an irreducible variety of dimension $\le\dim(G^0/L)+\dim S$.
Since $Y_{L,S}$ is an irreducible constructible subset of dimension 
$\dim(G^0/L)+\dim S$ of $Y'_{L,S}$ (see 3.13(b)), it follows that $Y_{L,S}$ is dense in
$Y'_{L,S}$. The lemma is proved.

\proclaim{Proposition 3.15}Let $(L,S)\in\AA$. The closure of $Y_{L,S}$ in $G$ is a 
union of strata of $G$.
\endproclaim
Let $P,\bS,Y'_{L,S}$ be as in 3.14. By Lemma 3.14, it is enough to show that $Y'_{L,S}$
is a union of strata of $G$. Since $Y'_{L,S}$ is stable under $G^0$-conjugacy, it is 
enough to show:

{\it if $g\in Y'_{L,S}$, $g'\in G$, $g\si g'$, then $g'\in Y'_{L,S}$}
\nl
or the stronger statement:
 
(a) {\it if $g\in Y'_{L,S}$ and $z\in T_G(g)$ then $zg\in Y'_{L,S}$.}
\nl
Replacing $(L,S,P)$ by a $G^0$-conjugate we may assume that $g\in\bS U_P$. Now 
$g_s\in N_GP$ is semisimple hence it normalizes some Levi of $P$, see 1.4(a). Hence, 
replacing $(L,S)$ by a $U_P$-conjugate we may assume in addition that 
$g_s\in N_GL\cap N_GP$. Let $f$ be the obvious projection of the semidirect product 
$(N_GL\cap N_GP)U_P$ (see 1.26) onto $N_GL\cap N_GP$ (a homomorphism of algebraic 
groups). Let $h=f(g)$. We have $g=hu$ where $h\in\bS,u\in U_P$ and $h_s=f(g_s)$. Since 
$g_s\in N_GL\cap N_GP$ we have $f(g_s)=g_s$ so that $g_s=h_s$. Since the set of 
isolated elements of $N_GL$ is closed in $N_GL$ (see Lemma 2.8) and it contains $S$, it
must also contain $\bS$; thus, $h$ is isolated in $N_GL$ so that

(b) $T_{N_GL}(h)=(\cz_L\cap Z_G(h))^0$.
\nl
We show that 

(c) $T_G(g)=T_G(h)$. 
\nl
Using 2.1(e) and the equality $g_s=h_s$ we see that it is enough to show that 
$u\in Z_G(g_s)^0$. Now $g_s=h_s$ commutes with $g$ and $h$ hence it commutes with $u$. 
Thus $u\in U_P\cap Z_G(g_s)=U_P\cap Z_G(g_s)^0$ (see 1.11) and (c) follows. We show 
that

(d) $T_G(h)\sub T_{N_GL}(h)$.
\nl
that is, $(\cz_{Z_G(h_s)^0}\cap Z_G(h_u))^0\sub(\cz_{Z_L(h_s)^0}\cap Z_G(h_u))^0$. It 
is enough to show that $\cz_{Z_G(h_s)^0}\sub\cz_{Z_L(h_s)^0}$. This follows from 
1.12(b) since $h_s\in N_GP$. (Since $S\sub N_GP$, we have $\bS\sub N_GP$ hence 
$h\in N_GP$.)

From (b),(c),(d) we deduce that $T_G(g)\sub(\cz_L\cap Z_G(h))^0$. Hence to prove (a) it
is enough to show:

{\it if $g=hu, h\in\bS,u\in U_P$ and $z\in(\cz_L\cap Z_G(h))^0$ then $zg\in\bS U_P$.}
\nl
It is enough to show that $z\bS\sub\bS$. This follows from $zS\sub S$ (see 2.3(a)). The
lemma is proved.

\proclaim{Proposition 3.16} For $(L,S)\in\AA$, $Y_{L,S}$ is a locally closed 
(irreducible) subvariety of $G$. In particular, $Y_{L,S}$ is open in $Y'_{L,S}$.
\endproclaim
This follows from Proposition 3.15 using the following general fact.

Assume that we are given an algebraic variety $V$ and a partition $V=\sqc_{j\in J}V_j$ 
where $V_j$ are irreducible constructible subsets of $V$ ($J$ is finite) such that for 
any $j\in J$, the closure of $V_j$ is a union of subsets of the form $V_{j'}$. Then 
each $V_j$ is locally closed in $V$.
\nl
We may assume that $J=\{1,2,\do,n\}$ and 
$j'\le j$ whenever $V_{j'}$ is contained in the closure of $V_j$. Then for any 
$j$, $\sqc_{j';j'\le j}V_{j'}$ and $\sqc_{j';j'<j}V_{j'}$ are closed in $V$ and 
$V_j=\sqc_{j';j'\le j}V_{j'}-\sqc_{j';j'<j}V_{j'}$ is locally closed in $V$.

\subhead 3.17\endsubhead
We show that, for $(L,S)\in\AA$, $Y_{L,S}$ {\it is a smooth variety.} Since $\tY_{L,S}$
is a principal $\cw_S$-bundle over $Y_{L,S}$ (see 3.13(a)), it is enough to show that 
$\tY_{L,S}$ is smooth. Consider the morphism $\tY_{L,S}@>>>G^0/L,(g,xL)\m xL$. This is 
a $G^0$-equivariant fibration over the homogeneous space $G^0/L$ whose fibre at $L$ is 
$S^*$, which is smooth (being open in the homogeneous space $S$). 

\head 4. Dimension estimates\endhead
\subhead 4.1\endsubhead
The estimates in this section are generalizations of results in \cite{\IC,\S1} which we
follow closely.

For any parabolic $P$ of $G^0$ we set $\tP=N_GP,\utP=\tP/U_P,\uP=P/U_P=\utP^0$; let 
$\p_P:\tP@>>>\utP$ be the canonical map. Let $\cp$ be a $G^0$-conjugacy class of 
parabolics of $G^0$. Assume that for each $P\in\cp$ we are given a $\uP$-conjugacy 
class $\boc_P$ in $\utP$ such that for any $P\in\cp$ and any $g\in G^0$, $\Ad(g)$ 
carries $\p_P\i(\boc_P)$ onto $\p_{gPg\i}(\boc_{gPg\i})$. Let
$$\zz=\{(g,P_1,P_2)\in G\T\cp\T\cp;
g\in\p_{P_1}\i(\cz_{\uP_1}^0\boc_{P_1})\cap\p_{P_2}\i(\cz_{\uP_2}^0\boc_{P_2})\},$$
$$\zz'=\{(g,P_1,P_2)\in G\T\cp\T\cp;g\in
\p_{P_1}\i(\boc_{P_1})\cap\p_{P_2}\i(\boc_{P_2})\}.$$
We have a partition $\zz=\cup_\co\zz_\co$ where $\co$ runs over the $G^0$-orbits in 
$\cp\T\cp$ and $\zz_\co=\{(g,P_1,P_2)\in\zz;(P_1,P_2)\in\co\}$. Similarly we have a
partition $\zz'=\cup_\co\zz'_\co$ where 
$\zz'_\co=\{(g,P_1,P_2)\in\zz';(P_1,P_2)\in\co\}$. We say that $\co$ is {\it good} if 
for some/any $(P_1,P_2)\in\co$, $P_1,P_2$ have a common Levi. We say that $\co$ is {\it
bad} if it is not good. Let $\nu$ be the number of positive roots of $G^0$. Let $\bnu$ 
be the number of positive roots of $\uP$, $\bc=\dim\boc_P$, 
$\br=\dim(\cz_{\uP}\cap Z_{\uP}(\g))^0$ for $P\in\cp,\g\in\boc_P$.

\proclaim{Proposition 4.2} Let $\boc$ be a $G^0$-conjugacy class in $G$, $c=\dim\boc$.

(a) For any $P\in\cp,x\in\boc_P$ we have $\dim(\boc\cap\p_P\i(x))\le(c-\bc)/2$.

(b) For any $g\in\boc$ we have 
$\dim\{P\in\cp;g\in\p_P\i(\boc_P)\}\le(\nu-\fra{c}{2})-(\bnu-\fra{\bc}{2})$.

(c) Let $d=2\nu-2\bnu+\bc+\br$. Then $\dim\zz_\co\le d$ if $\co$ is good and 
$\dim\zz_\co<d$ if $\co$ is bad. Hence $\dim\zz\le d$.

(d) Let $d'=2\nu-2\bnu+\bc$. Then $\dim\zz'_\co\le d'$ for all $\co$. Hence 
$\dim\zz'\le d$.
\endproclaim
(We make the convention that the empty set has dimension $-\iy$.) In the case where 
$\cp=\{G^0\}$, the proposition is trivial. Therefore we may assume that $\cp\ne\{G^0\}$
and that the result is already known when $G$ is replaced by a group of strictly
smaller dimension. 

We prove (c) and (d). We can map $\zz_\co$ and $\zz'_\co$ to $\co$ by 
$(g,P_1,P_2)\m(P_1,P_2)$. We see that proving (c) and (d) for $\zz_\co,\zz'_\co$ is the
same as proving that for a fixed $(P',P'')\in\co$ we have

(c${}'$) $\dim\{\p_{P'}\i(\cz_{\uP'}^0\boc_{P'})\cap
\p_{P''}\i(\cz_{\uP''}^0\boc_{P''})\}\le d-\dim\co$,

(d${}'$) $\dim\{\p_{P'}\i(\boc_{P'})\cap\p_{P''}\i(\boc_{P''})\}\le d'-\dim\co$,
\nl
with strict inequality in (c${}'$) if $\co$ is bad. Choose Levi subgroups $L'$ of $P'$ 
and $L''$ of $P''$ such that $L',L''$ contain a common maximal torus. Then $P'\cap P''$
is a connected group with Levi $L'\cap L''$. Let $\tL'=N_GL'\cap\tP'$,
$\tL''=N_GL''\cap\tP''$. By 1.26, we may identify $\tL'=\utP'$ via $\p_{P'}$. Similarly
we identify $\tL''=\utP''$ via $\p_{P''}$. Thus we regard 
$\boc_{P'}\sub\tL',\boc_{P''}\sub\tL''$. 
If $g\in\tP'\cap\tP''$ then, by 1.25(a),(b), we may 
write uniquely $g$ in the form $zu''u=zu'v$ where $z\in\tL'\cap\tL''$,
$u''\in L'\cap U_{P''},u\in U_{P'}\cap P'',u'\in L''\cap U_{P'},v\in U_{P''}\cap P'$.
We see that (c${}'$),(d${}'$) are equivalent to
$$\align &\dim\{(u,v,u'',u',z)\in(U_{P'}\cap P'')\T(U_{P''}\cap P')\T(L'\cap U_{P''})
\T(L''\cap U_{P'})\\& \T(\tL'\cap\tL'');u''u=u'v,zu''\in\cz_{L'}^0\boc_{P'},
zu'\in\cz_{L''}^0\boc_{P''}\}\le d-\dim\co,\tag c''\endalign$$
$$\align &\dim\{(u,v,u'',u',z)\in(U_{P'}\cap P'')\T(U_{P''}\cap P')\T(L'\cap U_{P''})
\T(L''\cap U_{P'})\\& \T(\tL'\cap\tL'');u''u=u'v,zu''\in\boc_{P'},zu'\in\boc_{P''}\}
\le d'-\dim\co,\tag d''\endalign$$
with strict inequality in (c${}''$) if $\co$ is bad. When 
$(u',u'')\in(L''\cap U_{P'})\T(L'\cap U_{P''})$ is fixed, the variety
$$R=\{(u,v)\in(U_{P'}\cap P'')\T(U_{P''}\cap P');u''u=u'v\}=
\{(u,v)\in U_{P'}\T U_{P''};u''u=u'v\}$$
is isomorphic to $U_{P'}\cap U_{P''}$. (Indeed, if we set
$\ti u=u'{}\i u''uu''{}\i\in U_{P'}$, $\ti v=vu''{}\i\in U_{P''}$, then $R$ becomes 
$\{(\ti u,\ti v)\in U_{P'}\T U_{P''};\ti u=\ti v\}$.) Since 
$\dim(U_{P'}\cap U_{P''})=2\nu-2\bnu-\dim\co$ we see that (c${}''$),(d${}''$) are 
equivalent to
$$\align &\dim\{(u'',u',z)\in(L'\cap U_{P''})\T(L''\cap U_{P'})\T(\tL'\cap\tL'');\\&
zu''\in\cz_{L'}^0\boc_{P'},zu'\in\cz_{L''}^0\boc_{P''}\} \le\bc+\br,\tag e
\endalign$$
$$\align &\dim\{(u'',u',z)\in(L'\cap U_{P''})\T(L''\cap U_{P'})\T(\tL'\cap\tL'');\\&
zu''\in\boc_{P'},zu'\in\boc_{P''}\}\le\bc,\tag f\endalign$$
with strict inequality in (e) if $\co$ is bad. 

Let us consider the variety in (f). Let $\p_3$ be the projection of that variety on the
$z$-coordinate. We show 

(g) $\text{image}(\p_3)$ {\it is a union of finitely many $L'\cap L''$-conjugacy classes 
in the reductive group $\tL'\cap\tL''$ with identity component} $L'\cap L''$.
\nl
Let $H'=\tL'\cap\tP''=(\tL'\cap\tL'')(L'\cap U_{P''})$ (semidirect product) and let
$f':H'@>>>\tL'\cap\tL''$ be the projection $zu''\m z$. Let 
$H''=\tL''\cap\tP'=(\tL'\cap\tL'')(L''\cap U_{P'})$ (semidirect product) and let
$f'':H''@>>>\tL''\cap\tL'$ be the projection $zu'\m z$. Then 
$\text{image}(\p_3)=f'(H'\cap\boc_{P'})\cap f''(H''\cap\boc_{P''})$.
Using 1.15(a) for $G$ or reductive groups of smaller dimension, we see that it is 
enough to show that 
$$\text{image}(\p_3)_s=f'(H'\cap\boc_{P'})_s\cap f''(H''\cap\boc_{P''})_s=
f'(H'\cap(\boc_{P'})_s)\cap f''(H''\cap(\boc_{P''})_s)$$
is a union of finitely many (semisimple) $L'\cap L''$-conjugacy classes in 
$\tL'\cap\tL''$. It is enough to show that $H'\cap(\boc_{P'})_s$ is 
a union of finitely 
many $H'{}^0=L'\cap P''$-conjugacy classes in $H'$ and that $H''\cap(\boc_{P''})_s$ is 
a union of finitely many $H''{}^0=L''\cap P'$-conjugacy classes in $H''$. Since 
$(\boc_{P'})_s$ is a semisimple $L'$-conjugacy class in $\tL'$ and $H'$ is a closed 
subgroup of $\tL'$, the intersection $H'\cap(\boc_{P'})_s$ is a union of finitely many 
$H'{}^0$-conjugacy classes in $H'$, see 1.27; similarly, $H''\cap(\boc_{P''})_s$ is a 
union of finitely many $H''{}^0$-conjugacy classes in $H''$. This proves (g). 

We can write $\text{image}(\p_3)=\c_1\cup\c_2\cup\do\cup\c_n$ where $\c_i$ are 
$(L'\cap L'')$-conjugacy classes in $\tL'\cap\tL''$. The inverse image under $\p_3$ of 
a point $z\in\c_i$ is a product of two varieties of the type considered in (a) but for
a smaller group ($G$ replaced by $\tL'$ or $\tL''$) hence by the induction hypothesis 
it has dimension $\le(\bc-\dim\c_i)/2+(\bc-\dim\c_i)/2$. Hence 
$\dim\p_3\i(\c_i)\le\bc$. Since this holds for each $i\in[1,n]$, we see that the
variety in (f) has dimension $\le\bc$. Thus, (d) is proved (assuming the induction
hypothesis).

We now consider the variety in (e). Let $\tp_3$ be the projection of that variety on 
the $z$-coordinate. With the earlier notation we have 
$$\text{image}(\tp_3)=f'(H'\cap\cz_{L'}^0\boc_{P'})\cap 
f''(H''\cap\cz_{L''}^0\boc_{P''}).$$
By 1.21(d) (for $\tL',\tL''$ instead of $G$) we have 
$$\cz_{L'}^0\boc_{P'}={}^{\d'}\cz_{L'}^0\boc_{P'},
\cz_{L''}^0\boc_{P''}={}^{\d''}\cz_{L''}^0\boc_{P'}$$
where $\d'$ (resp. $\d''$) is the connected component of $\tL'$ (resp. $\tL''$) that
contains $\boc_{P'}$ (resp. $\boc_{P''}$). We have 
${}^{\d'}\cz_{L'}^0\sub\cz_{L'}^0\sub\cz_{L'\cap L''}^0\sub H'$ and 
$f'(\z h)=\z f'(h)$ for $\z\in{}^{\d'}\cz_{L'}^0,h\in H'$ hence
$f'(H'\cap\cz_{L'}^0\boc_{P'})={}^{\d'}\cz_{L'}^0f'(H'\cap\boc_{P'})$. Similarly, 
$f''(H''\cap\cz_{L''}^0\boc_{P''})={}^{\d''}\cz_{L''}^0f''(H''\cap\boc_{P''})$. By an 
earlier argument we have $f'(H'\cap\boc_{P'})=\e'_1\cup\e'_2\cup\do\cup\e'_r$,
$f''(H''\cap\boc_{P''})=\e''_1\cup\e''_2\cup\do\cup\e''_t$ where 
$\e'_1,\e'_2,\do,\e'_r,\e''_1,\e''_2,\do,\e''_t$ are 
$(L'\cap L'')$-conjugacy classes in $\tL'\cap\tL''$. Thus,
$$\text{image}(\tp_3)=\cup_{i\in[1,r],j\in[1,t]}
({}^{\d'}\cz_{L'}^0\e'_i\cap({}^{\d''}\cz_{L''}^0\e''_j).$$
Here the set corresponding to $i,j$ is empty unless $\e'_i,e''_j$ are contained in the 
same connected component $X=X_{ij}$ of $\tL'\cap\tL''$. In that case we  have

${}^{\d'}\cz_{L'}^0\sub(\cz_{L'\cap L''}^X)^0$ and
${}^{\d''}\cz_{L''}^0\sub(\cz_{L'\cap L''}^X)^0$. 
\nl
(Indeed, since $L'\cap L''$ is a Levi of a parabolic of $L'$, we have 
$\cz_{L'}^0\cap\cz_{L'\cap L''}^0$. Let $x\in\e'_i\sub\tL'\cap\tL''$. We have 
$x=f'(\ti x)$ for some $\ti x\in H'\cap\cz_{L'}^0\boc_{P'}$. If 
$z\in{}^{\d'}\cz_{L'}^0$, then $z\in\cz_{L'\cap L''}^0$ and $z\ti x=\ti x$. Hence 
$f'(z)f'(\ti x)=f'(\ti x)f'(z)$ that is, $zx=xz$ so that $z\in\cz_{L'\cap L''}^X$. We 
see that ${}^{\d'}\cz_{L'}^0\sub\cz_{L'\cap L''}^X$ hence
${}^{\d'}\cz_{L'}^0\sub(\cz_{L'\cap L''}^X)^0$. Similarly,
${}^{\d''}\cz_{L''}^0\sub(\cz_{L'\cap L''}^X)^0$, as required.)

Using now 1.24(a) we see that 
${}^{\d'}\cz_{L'}^0\e'_i\cap{}^{\d''}\cz_{L''}^0\e''_j$ is a finite union of sets of 
the form $({}^{\d'}\cz_{L'}^0\cap{}^{\d''}\cz_{L''}^0)\d$ where $\d$ is a 
$(L'\cap L'')$-conjugacy class in $\tL'\cap\tL''$ contained in $X$. We see that
$$\text{image}(\tp_3)=\cup_{k\in[1,m]}({}^{\d'}\cz_{L'}^0\cap{}^{\d''}\cz_{L''}^0)\d_k
$$ 
where $\d_1,\d_2,\do,\d_m$ are $(L'\cap L'')$-conjugacy classes in 
$\tL'\cap\tL''$ such that for any $k$, the connected component $X_k$ of $\tL'\cap\tL''$
that contains $\d_k$ satisfies

${}^{\d'}\cz_{L'}^0\sub {}^{X_k}\cz_{L'\cap L''}^0$ and
${}^{\d''}\cz_{L''}^0\sub {}^{X_k}\cz_{L'\cap L''}^0$. 
\nl
The inverse image under $\tp_3$ of a point 
$z\in({}^{\d'}\cz_{L'}^0\cap{}^{\d''}\cz_{L''}^0)\d_k$ is a product of two varieties of
the type considered in (a) but for a smaller group ($G$ replaced by $\tL'$ or $\tL''$) 
hence by the induction hypothesis it has dimension 
$\le(\bc-\dim\d_k)/2+(\bc-\dim\d_k)/2$. Hence 
$$\align&\dim\tp_3\i(({}^{\d'}\cz_{L'}^0\cap{}^{\d''}\cz_{L''}^0)\d_k)
\le\bc -\dim\d_k+\dim({}^{\d'}\cz_{L'}^0\cap{}^{\d''}\cz_{L''}^0\d_k)\\&
=\bc -\dim\d_k+\dim({}^{\d'}\cz_{L'}^0\cap{}^{\d''}\cz_{L''}^0)+\dim\d_k\le
\bc+\dim({}^{\d'}\cz_{L'}^0)=\bc+\br\tag h\endalign$$
where the equality comes from 1.23(a). Since this holds for each $k\in[1,m]$, we see 
that the variety in (e) has dimension $\le\bc+\br$.

Moreover, if the second inequality in (h) is an equality then $\co$ is good. (Indeed, 
in this case we have ${}^{\d'}\cz_{L'}^0={}^{\d''}\cz_{L''}^0$. Taking centralizers in 
$G^0$ for both sides of the previous equality and using 1.10(a), we obtain $L'=L''$ 
hence $\co$ is good.) Thus (c) is proved (assuming the induction hypothesis).

We show that (b) is a consequence of (d). Let 
$\zz'(\boc)=\{(g,P_1,P_2)\in\zz';g\in\boc\}$. If $\zz'(\boc)=\em$ then the variety in 
(b) is empty and (b) follows. Hence we may assume that $\zz'(\boc)\ne\em$. From (d) we 
have $\dim\zz'(\boc)\le d'$. We map $\zz'(\boc)$ onto $\boc$ by the first projection. 
Each fibre of this map is a product of two copies of the variety in (b). It follows 
that the variety in (b) has dimension equal to
$$(\dim\zz'(\boc)-\dim\boc)/2\le(d'-c)/2=\nu-\bnu+\fra{\bc}{2}-\fra{c}{2}$$
and (b) is proved.

We show that (a) is a consequence of (b). Consider the variety
$\{(g,P)\in\boc\T\cp;g\in\p_P\i(\boc_P)\}$. Projecting it to the $g$-coordinate and 
using (b) we see that it has dimension $\le\nu-\bnu+\fra{\bc}{2}+\fra{c}{2}$. If we 
project it to the $P$-coordinate, each fibre will be isomorphic to the variety
$\boc\cap\p_P\i(\boc_P)$ (with $P\in\cp$ fixed). Hence
$$\dim(\boc\cap\p_P\i(\boc_P))\le\nu-\bnu+\fra{\bc}{2}+\fra{c}{2}-\dim\cp=
\fra{\bc}{2}+\fra{c}{2}.$$
Now $\boc\cap\p_P\i(\boc_P)$ maps onto $\boc_P$ (via $\p_P$) and each fibre is the 
variety in (a). Hence the variety in (a) has dimension 
$\le\fra{\bc}{2}+\fra{c}{2}-\bc=(c-\bc)/2$. The proposition is proved.

\subhead 4.3\endsubhead
In the case where $\co$ is good, the variety in 4.2(e) is
$\cz_{L'}^0\boc_{P'}\cap\cz_{L'}^0\boc_{P''}$ since $L'=L''$. This is empty if
$\cz_{L'}^0\boc_{P'}\ne\cz_{L'}^0\boc_{P''}$ and is $\cz_{L'}^0\boc_{P'}$ if 
$\cz_{L'}^0\boc_{P'}=\cz_{L'}^0\boc_{P''}$. In the last case, 
it follows that $\zz_\co$ is 
irreducible of dimension equal to $d$.

\subhead 4.4\endsubhead
The inequality in 4.2(b) can be reformulated as follows. Let us fix $P_0\in\cp$ and a
$\uP_0$-conjugacy class $\boc_0$ in $\utP_0$ with $\dim\boc_0=\bc$. Let $\d$ be the 
connected component of $\utP_0$ that contains $\boc_0$. Then for any $g\in\boc$ we have

(a) $\dim\{xP_0\in G^0/P_0; x\i gx\in\p_{P_0}\i(\boc_0)\}\le
(\nu-\fra{c}{2})-(\bnu-\fra{\bc}{2})$.
\nl
We have the following variant of (a):

(b) $\dim\{xP_0\in G^0/P_0;x\i gx\in\p_{P_0}\i(\cz_{\uP_0}^0\boc_0)\}\le
(\nu-\fra{c}{2})-(\bnu-\fra{\bc}{2})$.
\nl
This follows from (a) by observing that, for given $g$, there exist finitely many
$\uP_0$-conjugacy classes $\boc^1,\boc^2,\do,\boc^t$ in $\utP_0$ of dimension $\bc$ 
such that

$x\in G^0,x\i gx\in\p_{P_0}\i(\cz_{\uP_0}^0\boc_0)\imp 
x\i gx\in\p_{P_0}\i(\boc^1\cup\do\cup\boc^t)$.
\nl
Since $\cz_{\uP_0}^0\boc_0={}^\d\cz_{\uP_0}^0\boc_0$, it is enough to show that 

(c) ${}^\d\cz_{\uP_0}^0\boc_0\cap\p_{P_0}(\boc\cap\tP_0)$ 
\nl
is a union of finitely many $\uP_0$-conjugacy classes in $\utP_0$. (All 
$\uP_0$-conjugacy classes contained in ${}^\d\cz_{\uP_0}^0\boc_0$ have dimension $\bc$.)
Using 1.15(a) it is enough to show that the set of semisimple parts of the elements in 
(c) that is, ${}^\d\cz_{\uP_0}^0(\boc_0)_s\cap\p_{P_0}(\boc_s\cap\tP_0)$, is a finite 
union of (semisimple) $\uP_0$-conjugacy classes. This follow from the fact that 
$\boc_s\cap\tP_0$ is a finite union of (semisimple) $P_0$-conjugacy classes in $\tP_0$ 
(see 1.27).

\head 5. Some complexes on $G$\endhead
\subhead 5.1\endsubhead
We fix a prime number $l$ invertible in $\kk$. We use the term "local system" instead
of $\bbq$-local system. For an algebraic variety $V$ let $\cd(V)$ be the bounded
derived category of $\bbq$-sheaves on $V$. For $K\in\cd(V)$ let $\ch^iK$ be the
$i$-th cohomology sheaf of $K$.

\subhead 5.2\endsubhead
Let $C$ be an isolated stratum of $G$. For any integer $n\ge 1$, invertible in $\kk$, 
let $\cs_n(C)$ be the category whose objects are the local systems on $C$ that are
equivariant for the (transitive) $\cz_{G^0}^0\T G^0$-action 

(a) $(z,x):g@>>>xz^ngx\i$
\nl
on $C$.
\nl
If a local system is in $\cs_n(C)$ then it is also in $\cs_{n'}(C)$ for any $n'\ge 1$ 
invertible in $\kk$ such that $n'\in n\ZZ$. Let $\cs(C)$ be the category whose objects 
are the local systems on $C$ that are in $\cs_n(C)$ for some $n$ as above.

\proclaim{Lemma 5.3}Let $C$ be an isolated stratum of $G$. Let $C'$ be the image of $C$
under the obvious map $\p:G@>>>G_{ss}$ (an isolated stratum of $G_{ss}$). Let $C''$ be
the image of $C$ under the obvious map $\r:G@>>>G/G^0_{der}$, where $G^0_{der}$ is the 
derived group of $G^0$ ($C''$ is a connected component of $G/G^0_{der}$ and an isolated
stratum of $G/G^0_{der}$). The following two conditions for a local system $\ce$ on $C$
are equivalent:

(i) $\ce\in\cs(C)$;

(ii) $\ce\cong\op_{i=1}^m\p^*\ce'_i\ot\r^*\ce''_i$ where $\ce'_i\in\cs(C')$,
$\ce''_i\in\cs(C'')$ are irreducible and $\p:C@>>>C',\r:C@>>>C''$ are the restrictions 
of $\p,\r$ above.
\endproclaim
The fact that (ii) implies (i) is immediate. We prove the converse. Let ${}'G=G/\G$ 
where $\G=\cz_{G^0}^0\cap G^0_{der}$ (a finite normal subgroup of $G$). The image
${}'C$ of $C$ under the obvious map $\l:G@>>>{}'G$ is an isolated
stratum of ${}'G$ and $\l$ restricts to $\l:C@>>>{}'C$ (a principal covering with group
$\G$). Let $\ce$ be as in (i). Then $\l_*\ce\in\cs({}'C)$ (note however that, if 
$\ce\in\cs_n(C)$, then $\l_*\ce$ is not necessarily in $\cs_n({}'C)$). We may assume 
that $\ce$ is irreducible. There is natural surjective map $\l^*\l_*(\ce)@>>>\ce$. 
Since $\l^*\l_*(\ce)$ is semisimple, $\ce$ is a direct summand of $\l^*\l_*(\ce)$. We
can find an irreducible direct summand ${}'\ce$ of the local system $\l_*(\ce)$ on 
${}'C$ such that $\ce$ is a direct summand of $\l^*({}'\ce)$. We have canonically
${}'G^0=G^0/\cz_{G^0}^0\T G^0/G^0_{der}$ and ${}'C=C'\T C''$. It follows from the 
definitions that ${}'\ce$ is an external tensor product $\ce'\bxt\ce''$ where 
$\ce'\in\cs(C'),\ce''\in\cs(C'')$ are irreducible local systems.
Since $\ce''$ is irreducible and equivariant with respect to a transitive action of a 
connected commutative group, it has rank $1$. Since $\p^*\ce'$ is the inverse image of
an irreducible local system under a smooth map with connected fibres, it is itself an
irreducible local system. Since $\l^*({}'\ce)=\p^*\ce'\ot\r^*\ce''$ is a tensor product
of an irreducible local system and a rank $1$ local system, it is irreducible. It 
follows that $\ce=\l^*({}'\ce)=\p^*\ce'\ot\r^*\ce''$. The lemma is proved.

\subhead 5.4\endsubhead
Let $(L,S)\in\AA$. Let $P$ be a parabolic of $G^0$ with Levi $L$ such that 
$S\sub N_GP$. To simplify notation, in this section we set $Y=Y_{L,S},\tY=\tY_{L,S}$. 
As in 3.14, let $X=\{(g,xP)\in G\T G^0/P;x\i gx\in\bS U_P\}$; let $\ph:X@>>>G$ be the 
first projection. Let $f$ be the obvious projection of the semidirect product 
$(N_GL\cap N_GP)U_P$ (see 1.26) onto $N_GL\cap N_GP$ (a homomorphism of algebraic 
groups). 

\proclaim{Lemma 5.5} $(g,xL)\m(g,xP)$ is an isomorphism $\g:\tY@>\si>>\ph\i(Y)$.
\endproclaim
We verify this only at the level of sets. Assume that $(g,xL),(g',x'L)\in\tY$ have the 
same image under $\ph$. Then $g=g'$ and $x'=xp$ with $p\in P$. We have 
$x\i gx\in S^*,x'{}\i gx'\in S^*$ hence $p\i x\i gxp\in S^*$. It follows that 
$L(x\i gx)=L=L(p\i x\i gxp)=p\i L(x\i gx)p=p\i Lp$. Thus, $p\i Lp=L$ so that $p\in L$ 
and $xL=x'L$. Thus, $\g$ is injective. 

To show that $\g$ is surjective it is enough to show that, if $g\in S^*,x\in G^0$ 
satisfy $x\i gx\in\bS U_P$, then $u\i x\i gxu\in S^*$ for some $u\in U_P$ or 
equivalently that, 

if $g'\in\bS U_P,x\in G^0$ satisfy $xg'x\i\in S^*$, then $u\i g'u\in S^*$ for some 
$u\in U_P$. 
\nl
Now $g'_s\in N_GP$ is semisimple hence it normalizes some Levi of $P$ (see 1.4(a))
that is, some $U_P$-conjugate of $L$. Hence, replacing $g',x$ by $u'{}\i g'u',xu'$ for 
some $u'\in U_P$ we may assume in addition that $g'_s\in N_GL\cap N_GP$. We have 
$g'=hv$ where $h=f(g')\in\bS,v\in U_P$ and $h_s=f(g'_s)$. Since 
$g'_s\in N_GL\cap N_GP$ we have $f(g'_s)=g'_s$ so that $g'_s=h_s$. Then 
$h\i g'\in U_P\cap Z_G(g'_s)=U_P\cap Z_G(g'_s)^0$. Using 2.1(e), we see that 
$T(g')=T(h)$. By 1.22(b), we can find $h'\in S$ such that $h_s=h'_s$ and 
$h'{}\i h\in Z_G(h_s)^0$. Then $T(h)=T(h')$, by 2.1(e). Thus, $T(g')=T(h')$. By 3.8(b)
we have $L\sub L(h')=L(g')$. Since $xg'x\i\in S^*$, we have $L(xg'x\i)=L$ hence 
$L(g')=x\i Lx$. Thus, $L\sub x\i Lx$. Since $L,x\i Lx$ are irreducible of the same 
dimension, we have $L=x\i Lx$. Hence $N_GL=x\i N_GLx$. Since $xg'x\i\in S^*\sub N_GL$,
we have $g'\in N_GL$. Thus, $g'\in N_GL\cap\bS U_P$ hence $g'\in\bS$. Since 
$g'\in x\i S^*x$, we see that $x\i S x\cap\bS\ne\em$. Now $x\i S x$ is a stratum of 
$N_GL$ since $x\in N_GL$. Since $\bS$ is a union of strata of $N_GL$ one of which is 
$S$ and the others have dimension $<\dim S$, we see that $x\i S x=S$. This, together 
with $x\in N_GL$, implies that $x\i S^*x=S^*$. Since $g'\in x\i S^*x$, we see that 
$g'\in S^*$. Thus, $\g$ is surjective. The lemma is proved.

\subhead 5.6\endsubhead
Let $\ce\in\cs(S)$. We define a local system $\tce$ on $\tY$ by the requirement that 
$b^*\ce=a^*\tce$ where $a(g,x)=(g,xL),b(g,x)=x\i gx$ in the diagram
$$\tY@<a<<\{(g,x)\in G\T G^0;x\i gx\in S^*\}@>b>>S.$$
(We use the fact that $a$ is a principal $L$-bundle and $b^*\ce$ is $L$-equivariant.)
For any stratum $S'$ of $N_GL\cap N_GP$ such that $S'\sub\bS$ we set
$X_{S'}=\{(g,xP)\in G\T G^0/P;x\i gx\in S'U_P\}$. Then $X=\sqc_{S'}X_{S'}$ 
(union over all $S'\sub\bS$ as above; there are only finitely many $S'$ in the union,
see 3.7). By Lemma 2.8, each $S'$ is an isolated stratum of $N_GL\cap N_GP$. Note that
each $X_{S'}$ is smooth, irreducible. We define a local system $\bce$ on $X_S$ by the 
requirement that $b'{}^*\ce=a'{}^*\bce$ where $a'(g,x)=(g,xP)$, 
$b'(g,x)=f(x\i gx)$ in the diagram
$$X_S@<a'<<\{(g,x)\in G\T G^0;x\i gx\in S U_P\}@>b'>>S$$
(We use the fact that $a'$ is a principal $P$-bundle and $b'{}^*\ce$ is 
$P$-equivariant.) It is easy to see that the restriction of $\bce$ to $\tY$ (identified
with an open subset of $X$ as in Lemma 5.5) is $\tce$. The intersection cohomology 
complexes $IC(X,\bce)$ (on $X$) and $IC(\bS,\ce)$ (on $\bS$) are related by 

(a) $b''{}^*IC(X,\bce)=a''{}^*IC(\bS,\ce)$
\nl
where $a''(g,x)=(g,xP),b''(g,x)=f(x\i gx)$ in the diagram
$$X@<a''<<\{(g,x)\in G\T G^0;x\i gx\in\bS U_P\}@>b''>>\bS.$$
Here $a''$ is a principal $P$-bundle and $b''$ is a locally trivial fibration with 
smooth connected fibres. We write $\ph:X@>>>\bY$ for the restriction of $\ph:X@>>>G$.
Here $\bY$ is the closure of $Y$ in $G^1$. Recall that we have a finite covering 
(principal bundle) $\p:\tY@>>>Y$ (see 3.13(a)) hence $\p_!\tce$ is a well defined local
system on $Y$. Thus $IC(\bY,\p_!\tce)$ is well defined (on $\bY$).

\proclaim{Proposition 5.7} $\ph_!(IC(X,\bce))$ is canonically isomorphic to
$IC(\bY,\p_!\tce)$.
\endproclaim
Let $K=IC(X,\bce)$ and let $K^*=IC(X,\bce^*)$ where $\bce^*$ is defined like $\bce$ by
replacing $\ce$ by the dual local system $\ce^*$. Then $K^*$ is the Verdier dual of $K$
with a suitable shift. Since $\ph$ is proper, it follows that $\ph_!(K^*)$ is the 
Verdier dual of $\ph_!K$ with a suitable shift. We have $K|_{\tY}=\bce|_{\tY}=\tce$. 
Using Lemma 5.5, we see that $\ph_!K|_Y=\p_!\tce$. By the definition of an intersection
cohomology complex we see that it is enough to verify the following statement.

For any $i>0$ we have $\dim\supp\ch^i(\ph_!K)<\dim Y-i$ and
$\dim\supp\ch^i(\ph_!(K^*))<\dim Y-i$.
\nl
We shall only verify this for $K$; the corresponding statement for $K^*$ is entirely
analogous.

If $g\in\bY$, the stalk $\ch^i_g(\ph_!K)$ at $g$ is equal to $H^i_c(\ph\i(g),K)$. We 
have a partition $\ph\i(g)=\cup_{S'}(\ph\i(g)\cap X_{S'})$ where $S'$ runs over the 
strata of $N_GL\cap N_GP$ contained in $\bS$. If $H^i_c(\ph\i(g),K)\ne 0$ then 
$H^i_c(\ph\i(g)\cap X_{S'},K)\ne 0$ for some $S'$. Hence it is enough to prove:

For any $i>0$ and any $S'$ as above we have 
$\dim\{g\in\bY;H^i_c(\ph\i(g)\cap X_{S'},K)\ne 0\}<\dim Y-i$.
\nl
Assume first that $S'\ne S$. If $H^i_c(\ph\i(g)\cap X_{S'},K)\ne 0$ then the 
hypercohomology spectral sequence for $K$ on $\ph\i(g)\cap X_{S'}$  shows that we can
write $i=j_1+j_2$ with $j_2\le 2\dim(\ph\i(g)\cap X_{S'})$ and
$\ch^{j_1}(K|_{\ph\i(g)\cap X_{S'}})\ne 0$ hence $\ch^{j_1}(K)|_{X_{S'}}\ne 0$. Using 
5.6(a), we see that $\ch^{j_1}(K)$ is a local system on $X_{S'}$ (since
$\ch^{j_1}IC(\bS,\ce)$ is a local system on $S'$, which is a $(\cz_L^0\T L)$-orbit on
$N_GL\cap N_GP$). Thus, $X_{S'}\sub\supp\ch^{j_1}(K)$. Since $K=IC(X,\bce)$, it follows
that $j_1<\dim X-\dim X_{S'}=\dim S-\dim S'$. Thus we have
$i<2\dim(\ph\i(g)\cap X_{S'})+\dim S-\dim S'$ and it is enough to show that
$$\dim\{g\in\bY;\dim(\ph\i(g)\cap X_{S'})
>\fra{i}{2}-\fra{1}{2}(\dim S-\dim S')\}<\dim Y-i.$$
If this is violated for some $i>0$, it would follow that the space of triples
$$\{(g,xP,x'P)\in\bY\T G^0/P\T G^0/P;x\i gx\in S'U_P,x'{}\i gx'\in S'U_P\}$$
has dimension $>\dim Y-i+i-(\dim S-\dim S')=\dim G^0/L+\dim S'$. This contradicts 
4.2(c).

Next, assume that $S'=S$. If $H^i_c(\ph\i(g)\cap X_S,K)\ne 0$ then 
$i\le 2\dim(\ph\i(g)\cap X_S)$ since $K|_{\ph\i(g)\cap X_S}$ is a local system. Hence
it is enough to show that for $i>0$ we have
$$\dim\{g\in\bY;\dim(\ph\i(g)\cap X_S)\ge\fra{i}{2}\}<\dim Y-i.$$
Assume that this is violated for some $i>0$. Thus, setting
$F=\{g\in\bY;\dim(\ph\i(g)\cap X_S)\ge\fra{i}{2}\}$, we have $\dim F\ge\dim Y-i$. Then
the space of triples 
$$\{(g,xP,x'P)\in F\T G^0/P\T G^0/P;x\i gx\in SU_P,x'{}\i gx'\in SU_P\}$$
has dimension $\ge\dim G^0/L+\dim S$ (and the last inequality is strict if
$\dim F>\dim Y-i$). From 4.2(c) we see that this space of triples has dimension 
$\le\dim G^0/L+\dim S$ hence it has dimension equal to $\dim G^0/L+\dim S$, which 
forces $\dim F=\dim Y-i$. We partition our space of triples into subsets by specifying 
the $G^0$-orbit of $(xPx\i,x'Px'{}\i)$ (for simultaneous conjugation).
By 4.2(c), the subset corresponding to a bad orbit has dimension $<\dim G^0/L+\dim S$. 
It follows that the subset corresponding to some good orbit has dimension equal to 
$\dim G^0/L+\dim S$. Thus there exists $n\in N_{G^0}L$ such that 
$$\{(g,xP,x'P)\in F\T G^0/P\T G^0/P;x\i gx\in SU_P,x'{}\i gx'\in SU_P;x\i x'\in PnP\}
\tag a$$
has dimension equal to $\dim G^0/L+\dim S$. By 4.3, the variety
$$\{(g,xP,x'P)\in G\T G^0/P\T G^0/P;x\i gx\in SU_P,x'{}\i gx'\in SU_P;x\i x'\in PnP\}
\tag b$$
is empty if $nSn\i\ne S$ and is irreducible of dimension $\dim G^0/L+\dim S$ if 
$nSn\i=S$. It follows that we must have $nSn\i=S$ and the variety (a) is dense in the
variety (b). It follows that $F$ is dense in the image $I$ of the variety (b) under the
first projection. Hence $\dim I=\dim Y-i$. If $g\in S^*$ then $(g,P,nP)$ belongs to the
variety (b). Thus, $S^*\sub I$. Since $I$ is stable under $G^0$-conjugacy, we must have
$Y\sub I$. Thus, $\dim I\ge\dim Y$. It follows that $\dim Y-i\ge\dim Y$ hence $i\le 0$,
contradicting $i>0$. The proposition is proved.

\head 6. Cuspidal local systems\endhead
\subhead 6.1\endsubhead
In this section we fix an isolated stratum $C$ of $G$. Let $D$ be the connected
component of $G$ that contains $C$.

\proclaim{Lemma 6.2} Let $P$ be a parabolic of $G$ and let $g\in C\cap N_GP$. Let 
$\boc_P$ be the $P/U_P$-conjugacy class of the image of $g$ in $N_GP/U_P$. Let
$\d=\dim C-\dim{}^D\cz_{G^0}^0-\dim\boc_P$. Then $\dim(C\cap gU_P)\le\d/2$. Hence for 
any $\ce\in\cs(C)$ we have $H^i_c(C\cap gU_P,\ce)=0$ for $i>\d$.
\endproclaim
If $y\in gU_P$, then the semisimple elements $y_s,g_s$ normalize $U_P$ and are in the 
same $U_P$-coset hence, by a standard argument, are $U_P$-conjugate. Using the 
finiteness of the number of unipotent classes in $Z_G(g_s)$, see 1.15, we deduce that 
$gU_P$ is contained in the union of finitely many $G^0$-conjugacy classes in $G$. It is
then enough to show that for any $G^0$-conjugacy class $\boc$ in $G$ such that 
$\boc\sub C$, we have $\dim(\boc\cap gU_P)\le\fra{1}{2}\d$. This follows from 4.2(a), 
since $\dim C=\dim\boc+\dim{}^D\cz_{G^0}^0$, see 1.23(b).

\subhead 6.3\endsubhead
Let $\ce\in\cs(C)$. We say that $\ce$ is a {\it cuspidal local system} or that 
$(C,\ce)$ is a {\it cuspidal pair} for $G$ if, for any $P,g$ as in 6.2, with 
$P\ne G^0$ we have $H^\d_c(C\cap gU_P,\ce)=0$ where $\d$ is as in 6.2.

\proclaim{Lemma 6.4} Let $\ce\in\cs(C)$. Let us write
$\ce=\op_{i=1}^m\p^*\ce'_i\ot\r^*\ce''_i$ as in Lemma 5.3. Then $(C,\ce)$ is a cuspidal
pair for $G$ if and only if $(C',\ce'_i)$ is a cuspidal pair for $G_{ss}$ for 
$i=1,\do,m$ (with $C'$ as in Lemma 5.3.)
\endproclaim
Since the property of being cuspidal is preserved by taking direct sums of local 
systems or by passage to a direct summand, we see that we may assume that $m=1$. Next 
we note that, for $P,g$ as in 6.2, we have
$(\p^*\ce'_1\ot\r^*\ce''_1)|_{C\cap gU_P}\cong(\p^*\ce'_1)|_{C\cap gU_P}$ since 
$(\r^*\ce''_1)|_{C\cap gU_P}\cong\bbq$. Hence we may also assume that $\ce''_1=\bbq$. 
We set $\ce'=\ce'_1$. We must show that $(C,\ce)$ is a cuspidal pair for $G$ if and 
only if $(C',\ce')$ is a cuspidal pair for $G_{ss}$.

Assume that $(C',\ce')$ is a cuspidal pair for $G_{ss}$. To show that $(C,\ce)$ is a 
cuspidal pair for $G$ we must show that for any $P,g$ as in 6.2 with $P\ne G^0$ we have
$H^\d_c(C\cap gU_P,\ce)=0$ where $\d$ is as in 6.2. Let $\bP=\p(P)$, a proper parabolic
of $G^0_{ss}$. By assumption we have $H^{\ti\d}_c(C'\cap\p(g)U_{\bP},\ce')=0$ where 
$C_{\bP}$ is the $\bP/U_{\bP}$-conjugacy class of the image of $\p(g)$ in 
$N_{G_{ss}}\bP/U_{\bP}$ and $\ti\d=\dim C'-\dim C_{\bP}$. It is clear that $\p$
restricts to an isomorphism $C\cap gU_P@>\si>>C'\cap\p(g)U_{\bP}$. Since 
$\dim C-\dim{}^D\cz_{G^0}^0=\dim C'$ and $\dim C_P=\dim C_{\bP}$ ($C_P$ as in 6.2), we 
have $\d=\ti\d$. Hence 
$H^\d_c(C\cap gU_P,\ce)\cong H^{\ti\d}_c(C'\cap\p(g)U_{\bP},\ce')$ and we see that 
$(C,\ce)$ is a cuspidal pair for $G$. The reverse implication is proved in a similar 
way. The lemma is proved.

\subhead 6.5\endsubhead
Lemma 6.4 shows that the study of ccuspidal pairs for $G$ can be reduced to
the analogous problem for $G_{ss}$. We will show that we can further reduce to the case
of a unipotent class.

Assume that $G=G_{ss}$. Let $\ce\in\cs(C)$. For any $x\in C_s$ let 
$\cc^x=\{u\in Z_G(x);u\text{ unipotent, }xu\in C\}$. Then $Z_{G^0}(x)$ acts 
transitively (by conjugation) on $\cc^x$. Let $\ce^x$ be the inverse image of $\ce$ 
under $\cc^x@>>>C,u\m xu$.

\proclaim{Lemma 6.6} The following three conditions are equivalent:

(i) $(C,\ce)$ is a cuspidal pair for $G$;

(ii) there exists $x\in C_s$ such that for some/any $Z_G(x)^0$-conjugacy class $\cc'$ 
in $\cc^x$, the pair $(\cc',\ce^x|_{\cc'})$ is cuspidal for $Z_G(x)$.

(iii) for any $x\in C_s$ and for some/any $Z_G(x)^0$-conjugacy class $\cc'$ in $\cc^x$,
the pair $(\cc',\ce^x|_{\cc'})$ is cuspidal for $Z_G(x)$.
\endproclaim
We prove (i) assuming that (ii) holds. Assume that we are given $g\in C$ and a proper 
parabolic $P$ of $G^0$ such that $g\in N_GP$. The $P/U_P$-conjugacy class of the image 
$\bg$ of $g$ in $N_GP/U_P$ is denoted by $\boc_P$. We must show that 
$H^\d_c(C\cap gU_P,\ce)=0$ where $\d=\dim C-\dim\boc_P$. Since $g\in N_G(U_P)$, we have
$g_s\in N_G(U_P)$ hence $(g_sU_P)_s=\cv$ where $\cv=\{vg_sv\i;v\in U_P\}\sub C_s$. 
Hence $y\m y_s$ is a morphism $f:C\cap gU_P@>>>\cv$. Since $f$ commutes with the 
conjugation action of $U_P$ on $C\cap gU_P$ and $\cv$, we have
$H^\d_c(C\cap gU_P,\ce)\cong H^{\d-2\dim\cv}_c(f\i(g_s),\ce)$. Now $u\m g_su$ defines 
$\{u\in\cc^{g_s};u\in g_uU_P\}@>\si>>f\i(g_s)$ or equivalently,
$\cc^{g_s}\cap g_uU_Q@>\si>>f\i(g_s)$ where $Q=P\cap Z_G(g_s)^0$ (a parabolic of 
$Z_G(g_s)^0$ with $U_Q=U_P\cap Z_G(g_s)=U_P\cap Z_G(g_s)^0$). It is then enough to show
that $H^{\d-2\dim\cv}_c(\cc^{g_s}\cap g_uU_Q,\ce')=0$. (We have 
$g_u\in N_GQ\cap Z_G(g_s)$.) Since $Z_{G^0}(g_s)$ acts transitively on $\cc^{g_s}$, we 
see that $\cc^{g_s}$ is a disjoint union of finitely many $Z_G(g_s)^0$-conjugacy 
classes. It is enough to show that for any such conjugacy class $\cc'$ we have 
$H^{\d-2\dim\cv}_c(\cc'\cap g_uU_Q,\ce')=0$. (If this holds for some $\cc'$ then it 
automatically holds for all $\cc'$, by the transitivity of the $Z_{G^0}(g_s)$-action on
$\cc^{g_s}$.) This would follow from (iii) provided that we verify:

$Q\ne Z_G(g_s)^0$ and $\d-2\dim\cv=\d'$
\nl
where $\boc_Q$ is the $Q/U_Q$-conjugacy class of the image of $g_u$ in 
$N_{Z_G(g_s)}Q/U_Q$ and $\d'=\dim\cc'-\dim\boc_Q$. 

If $Q=Z_G(g_s)^0$ then $Z_G(g_s)^0\sub P$. Hence $g$ is not isolated, a contradiction.

We now show that $\d-2\dim\cv=\d'$ that is, 
$\dim C-\dim\boc_P-2\dim\cv=\dim\cc'-\dim\boc_Q$. Now $Q/U_Q=Z_{P/U_P}(\bg_s)^0$ and 
$\boc_Q$ may be identified with the $Z_{P/U_P}(\bg_s)^0$-conjugacy class of $\bg_u$ in 
$N_GP/U_P$. Consider the morphism $x\m x_s,\boc_P@>>>(\boc_P)_s$; the fibre of this 
morphism at $\bg_s$ is the $Z_{P/U_P}(\bg_s)$-conjugacy class of $\bg_u$ in $N_GP/U_P$ 
which has pure dimension $\dim\boc_Q$. We see that 
$\dim\boc_P-\dim\boc_Q=\dim(\boc_P)_s$. We also have $\dim C=\dim\cc'+\dim C_s$. Thus
$\dim C-\dim\boc_P-2\dim\cv-\dim\cc'+\dim\boc_Q=\dim C_s-\dim(\boc_P)_s-2\dim\cv$.
To prove that this is $0$ it is enough to show that
$(\dim G-\dim Z_G(g_s))-(\dim L_1-\dim Z_{L_1}(g_s))-2(\dim U_P-\dim Z_{U_P}(g_s))=0$
where $L_1$ be a Levi subgroup of $P$ normalized by $g_s$. Since 
$\dim G=\dim L_1+2\dim U_P$ it is enough to show that
$\dim Z_G(g_s)=Z_{L_1}(g_s)+2\dim Z_{U_P}(g_s)$. This follows from the fact that 
$Z_G(g_s)^0$ is reductive and $Z_{L_1}(g_s)^0$ is a Levi of a parabolic of $Z_G(g_s)^0$
with unipotent radical $Z_{U_P}(g_s)$ (see 1.12(a)). This proves (i).

We prove (iii) assuming that (i) holds. Let $x\in C_s$ and let $\cc'$ be a 
$Z_G(x)^0$-conjugacy class in $\cc^x$. Assume that we are given $y\in\cc'$ and a proper
parabolic $Q$ of $Z_G(x)^0$ such that $y\in N_{Z_G(x)}Q$. The $Q/U_Q$-conjugacy class 
of the image $\bar y$ of $y$ in $N_{Z_G(x)}Q/U_Q$ is denoted by $\boc_Q$. We must show 
that $H^{\d'}_c(\cc'\cap yU_Q,\ce')=0$ where $\d'=\dim\cc'-\dim\boc_Q$. It is enough to
show that $H^{\d'}_c(\cc^{g_s}\cap g_uU_Q,\ce')=0$. Let $g=xy=yx$. We have $g\in C$ and
$x=g_s,y=g_u$. By 1.18(a),we can find a parabolic $P$ of $G^0$ such that $g\in N_GP$ 
and $P\cap Z_G(g_s)^0=Q$. Clearly, $P\ne G^0$. By the arguments and notation in the 
first part of the proof, we see that $H^{\d'}_c(\cc^{g_s}\cap g_uU_Q,\ce')$ is 
isomorphic to $H^\d_c(C\cap gU_P,\ce)$ which is $0$ by assumption. This proves (iii).

Now (ii) and (iii) are equivalent since $G^0$ acts transitively by conjugation on 
$C_s$. The lemma is proved.

\subhead 6.7\endsubhead
Let $A$ be a simple perverse sheaf on $G$. We say that $A$ is {\it admissible} if
the following condition is satisfied: there exists $(L,S)\in\AA$, a cuspidal
irreducible local system $\ce\in\cs(S)$ and an irreducible direct summand $\tce_1$
of the local system $\tce$ on $\bY$ (with $\bY,\tce$ defined as in 5.6 in terms of 
$L,S,\ce$), such that $A$ is isomorhic to $IC(\bY,\tce_1)[\dim\bY]$ regarded as a 
simple perverse sheaf on $G$ ($0$ on $G-\bY$).


\Refs
\widestnumber\key{L1}
\ref\key{\BO}\by A.Borel\book Linear algebraic groups\publ Benjamin\publaddr New York
\yr 1969\endref
\ref\key{\DS}\by J. de Siebenthal\paper Sur les groupes de Lie compactes non connexes
\jour Comment. Math. Helv.\vol 31\yr 1956 \pages 41-89\endref
\ref\key{\FU}\by G.Lusztig\paper On the finiteness of the number of unipotent classes
\jour Invent. Math.\vol 34\yr 1976\pages 201-213\endref
\ref\key{\IC}\by G.Lusztig\paper Intersection cohomology complexes on a reductive group
\jour Invent. Math.\vol 75\yr 1984\pages 205-272\endref
\ref\key{\CS}\by G.Lusztig\paper Character sheaves,I\jour Adv.Math.\vol 56\yr 1985
\pages 193-237\endref
\ref\key{\IN}\by G.Lusztig\paper Introduction to character sheaves\jour
Proc.Symp.Pure Math.\vol 47(1)\yr 1987 \linebreak
\publ Amer.Math.Soc.\pages 165-180\endref
\ref\key{\CL}\by G.Lusztig\paper Classification of unipotent representations of simple 
$p$-adic groups,II\jour Represent.Th.\vol 6\yr 2002\pages 243-289\endref
\ref\key{\SP} \by N.Spaltenstein\book Classes unipotents et sous-groupes de Borel,
Lecture Notes in Mathematics\publ Springer Verlag\publaddr New York\vol 946\yr 1982
\endref
\ref\key{\ST}\by R.Steinberg\book Endomorphisms of linear algebraic groups,
Memoirs of Amer.Math.Soc\vol 80\yr 1968\endref
\endRefs
\enddocument